\documentclass[preprint,showpacs,english,aps]{revtex4-1}
\usepackage{babel,rotating,dcolumn}
\usepackage{colordvi,graphicx,color,amsbsy,amsmath,bm,amssymb}

\begin{document}

\title{Volume penalization for inhomogeneous Neumann boundary conditions modeling scalar flux in complicated geometry 
}%
\author{Teluo Sakurai$^1$}
\author{Katsunori Yoshimatsu$^{2,}$\footnote{yoshimatsu@nagoya-u.jp}} 
\author{Naoya Okamoto$^{3,}$\footnote{Present address: Aichi Institute of Technology, 1247, Yachikusa, Yakusacho, Toyota, 470-0392, Japan}}
\author{Kai Schneider$^4$}

\affiliation{$^1$Department of Computational Science and Engineering, Nagoya University, Nagoya, 464-8603, Japan}
\affiliation{$^2$Institute of Materials and Systems for Sustainability, Nagoya University, Nagoya, 464-8603, Japan}
\affiliation{$^3$Center for Computational Science, Nagoya University, Nagoya, 464-8603, Japan}
\affiliation{$^4$Institut de Math\'ematiques de Marseille (I2M), Aix-Marseille Universit\'e, CNRS and Centrale Marseille, 39 rue F. Joliot-Curie, 13453 Marseille Cedex 13, France}


\date{Received 4 October 2018, Accepted 4 April 2019. { J. Comput. Phys. {\bf 390} 452--469 (2019)}
}
\begin{abstract}
We develop a volume penalization method for inhomogeneous Neumann boundary conditions, 
generalizing the flux-based volume penalization method for homogeneous Neumann boundary condition proposed by Kadoch et al. [J. Comput. Phys. 231 (2012) 4365].
The generalized method allows us to model scalar flux through walls in geometries of complex shape using simple, e.g. Cartesian, domains for solving the governing equations. 
We examine the properties of the method,
by considering a one-dimensional Poisson equation with different Neumann boundary conditions. 
The penalized Laplace operator is discretized by second order central finite-differences and interpolation.
The discretization and penalization errors are thus assessed for several test problems.
Convergence properties of the discretized operator and the solution of the penalized equation are analyzed.
The generalized method is then applied to an advection-diffusion equation coupled with the Navier--Stokes equations in an annular domain which is immersed in a square domain.
The application is verified by numerical simulation of steady free convection in a concentric annulus heated through the inner cylinder surface using an extended square domain.
\end{abstract}

\maketitle


\section{Introduction} \label{Sec1}
Numerical modeling of multiphysics problems in complicated geometries is still a challenging problem in computational fluid dynamics. 
For example, the numerical prediction of the exchanges of scalar quantities, e.g., heat and mass, between solid and fluid phases has numerous industrial applications. 
A non-exhaustive list includes radiators and heat exchangers, for example the fan cooler of CPUs, or crystal growth processes with important applications for semi-conductors and light emission devices. 
The heat or mass flux through interfaces between solid and fluid phases can be mathematically modeled, 
and leads typically to Neumann or Robin boundary conditions of the governing advection-diffusion equations.

Immersed boundary methods yield an attractive approach for solving partial differential equations in domains of complex shape. 
The underlying idea is to consider a simple computational domain for which effective numerical methods are available and the boundary conditions are imposed by so-called penalty or direct forcing terms; 
for a review we refer to Refs. \cite{Peskin,Mittal,reviewKai}.  
For Dirichlet boundary conditions, there is an abundant literature available. 
Penalization techniques for Neumann and Robin boundary conditions are more recent and less developed,
see e.g., Refs. \cite{KaKoAnSc,KoNgSc,Morales,Vasilyev, Ren2013,BCL2015}.
Most immersed boundary methods result in low order, i.e., first or second order, approximation and computational stiffness.

In this paper, we focus on the volume penalization (VP) method. 
In the pioneering work by Angot et al. \cite{Angot1999}, 
the VP method was developed for imposing Dirichlet boundary conditions on the Navier--Stokes equations.
The boundary conditions for velocity of viscous flow are given by no-slip conditions on the surface of the solid, 
such as walls and obstacles. 
The VP method models the solid as a porous medium whose permeability $\eta \, (>0)$ is sufficiently small.
Mathematically, it was shown that the solution of the penalized Navier--Stokes equations converges towards the solution of the Navier--Stokes equations with no-slip boundary conditions, 
as $\eta \rightarrow 0$ \cite{Angot1999, Carbou}.
This VP method has been applied to various flows, e.g., confined hydrodynamic turbulence \cite{KaiMarie}, 
confined magnetohydrodynamic turbulence \cite{Morales,KaiNeffaa}, fluid-structure interaction for moving obstacles \cite{Kolmenskiy} and for flexible beams \cite{KoEnSc}, and the aerodynamics of insect flight \cite{KoMoFaSc}.

A VP method for general boundary conditions of Neumann and Robin types was proposed by generalization of the VP method using a weak formulation \cite{Ramiere1}, 
and applied in the context of finite element methods and finite volume methods \cite{Ramiere2}. 
Kadoch et al. \cite{KaKoAnSc} extended this approach for pseudo-spectral discretizations using a strong formulation in the case of homogeneous Neumann boundary conditions.

The aim of the current work is to extend the VP method, developed in Ref. \cite{KaKoAnSc} for homogeneous Neumann boundary conditions to inhomogeneous ones. 
Thus net scalar flux can be modeled keeping the continuity of its flux. 
The proposed method is analyzed in detail for a one-dimensional (1D) Poisson equation where exact solutions for both the penalized equation and the non-penalized equation are available. 
The penalization error can be hence computed analytically. 
For the numerical discretization, second order finite-differences and interpolation are used. 
The discretization error is determined and a guide for choosing the spatial discretization and $\eta$ in actual implementations is presented. 

This paper is organized as follows.
In Section \ref{Sec2}, we develop a flux-based VP representation of inhomogeneous Neumann boundary conditions. 
In Section \ref{sec3}, the VP representation is extended to a two-dimensional (2D) penalized Poisson equation.
In Section \ref{sec5}, we apply the VP representation to an advection-diffusion equation coupled with Navier--Stokes equations in an annular domain.
Numerical results are given for 2D steady incompressible convection in an annulus subjected to heat flux through its inner wall.
Conclusions are given in Section \ref{sec6}.

\section{Poisson equation with inhomogeneous Neumann boundary conditions} \label{Sec2}
We consider a 1D Poisson equation with inhomogeneous Neumann boundary conditions modeled with the VP method. 
In this case, the exact solution of the non-penalized and the penalized equations can be obtained 
and thus the penalization error can be assessed. 
The discretization error of the penalized equation using finite differences is also analyzed.
\subsection{A flux-based volume penalization representation of inhomogeneous Neumann boundary conditions} \label{sec2-1}
We first consider a 1D Poisson equation,
\begin{equation}
-d_x^2 w(x) = f(x) \,\,\, {\mbox{for}} \,\,\, x \in \Omega_f \label{eq2-0}
\end{equation}
with a source term $f(x)$, which is here given by
\begin{equation}
f(x) = m^2 \cos{m x}. \label{eq2-1}
\end{equation}
Here, $m=1,2,\cdots$, and $d_x=d/dx$.
The fluid domain is $\Omega_f = \{x \ | \ 0 <x < \pi \}$ and inhomogeneous Neumann boundary conditions,
\begin{equation}
d_x w(0)=d_x w(\pi)=\alpha, \label{eq2-2}
\end{equation}
are imposed, where $\alpha $ is a real valued constant. 
The argument $x$ will be omitted, unless otherwise stated.

The compatibility condition, $ -[d_x w]_0^\pi = \int_0^\pi f (x) dx \, ( =0)  $, is satisfied. 
Therefore, exact solutions of Eq. (\ref{eq2-0}) with Eq. (\ref{eq2-1}) exist satisfying the condition (\ref{eq2-2}).
For simplicity, 
we have chosen here the same value $\alpha$ on the left and right boundaries, $x=0$ and $x=\pi$.
If $\alpha=0$, the boundary conditions (\ref{eq2-2}) are homogeneous.
In Section \ref{sec2-3} we will generalize this case to inhomogeneous Neumann boundary conditions with different values on the boundaries.
The exact solutions of Eqs. (\ref{eq2-0})--(\ref{eq2-2}) are given by
\begin{equation}
w(x) = \cos{m x}+\alpha x+A_0, \label{eq2-3}
\end{equation}
where $A_0$ is an additive constant which reflects the non-uniqueness of the solution.
To determine $A_0$, we set the condition $\int_{\Omega_f} w (x) d x=0$,
which yields the value $A_0=-\pi\alpha/2$.

Generalizing the VP representation for homogeneous Neumann boundary conditions \cite{KaKoAnSc},
we obtain the VP representation of Eqs. (\ref{eq2-0}) and (\ref{eq2-2}):
\begin{equation}
-d_x ( \theta d_x v + \chi \alpha)= (1-\chi) f, \label{eq2-4}
\end{equation}
where 
\begin{equation}
\theta (x) = 1-\chi (x) + \eta \chi (x),  \label{eq2-5}
\end{equation}
and $\eta \, (>0)$ is the penalization parameter.
The mask function $\chi (x)$ determines the geometry of $\Omega_f$, 
and is given by
\begin{eqnarray}
\chi(x)= \left \{
\begin{array}{lll}
0   & {\rm for} & 0<x<\pi ,  \\
1/2 & {\rm at}  & x=0, \pi , \\
1   & {\rm for} & \pi<x<2\pi .
\end{array}
\right. \label{eq2-7}
\end{eqnarray}
Here,
we have chosen to imbed the fluid domain $\Omega_f=\{ x \ | \ 0 < x < \pi \}$ into the larger domain $\Omega =\{x \ | \ 0 \le x < 2 \pi \} = {\bar \Omega_f} \cup \Omega_s$ with the solid domain $\Omega_s = \{ x \ | \ \pi <x < 2 \pi \}$, 
and ${\bar \Omega_f}$ denoting the closure of $\Omega_f$.
Then we can impose $2\pi$-periodic boundary conditions.
This VP representation (\ref{eq2-5}) is based on the continuity of the flux through the non-penalized boundaries at $x=0 \, (=2\pi)$ and $x=\pi$.
The flux at the fluid-solid interfaces satisfies
\begin{equation}
d_x v = \eta d_x v + \alpha \quad {\mathrm{at}} \quad x=0 \, (=2\pi),\pi . \label{eq2-8}
\end{equation}

\begin{figure}[tb]
\begin{center}
\includegraphics[width=7cm]{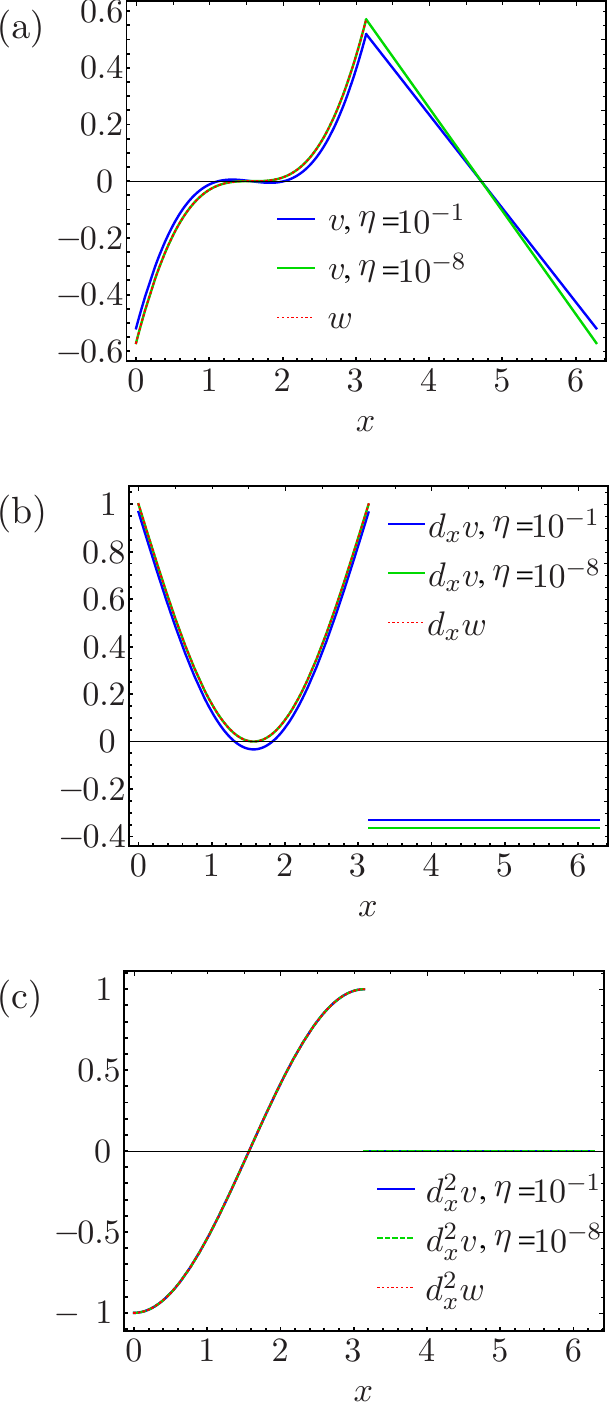} 
\end{center}
\caption{Comparison of the analytical solution $v(x)$, Eq. (\ref{eq2-9}), with the exact solution $w(x)$, Eq. (\ref{eq2-3}), for $\alpha=1$ at $\eta=10^{-1}$ and $10^{-8}$:
(a) $v$ and $w$, (b) the first derivatives $d_x v$ and $d_x w$, 
and (c) the second derivatives $d_x^2 v$ and $d_x^2 w$.
}\label{fig2-1}
\end{figure}

\begin{figure}[tb]
\begin{center}
\includegraphics[width=7cm,keepaspectratio]{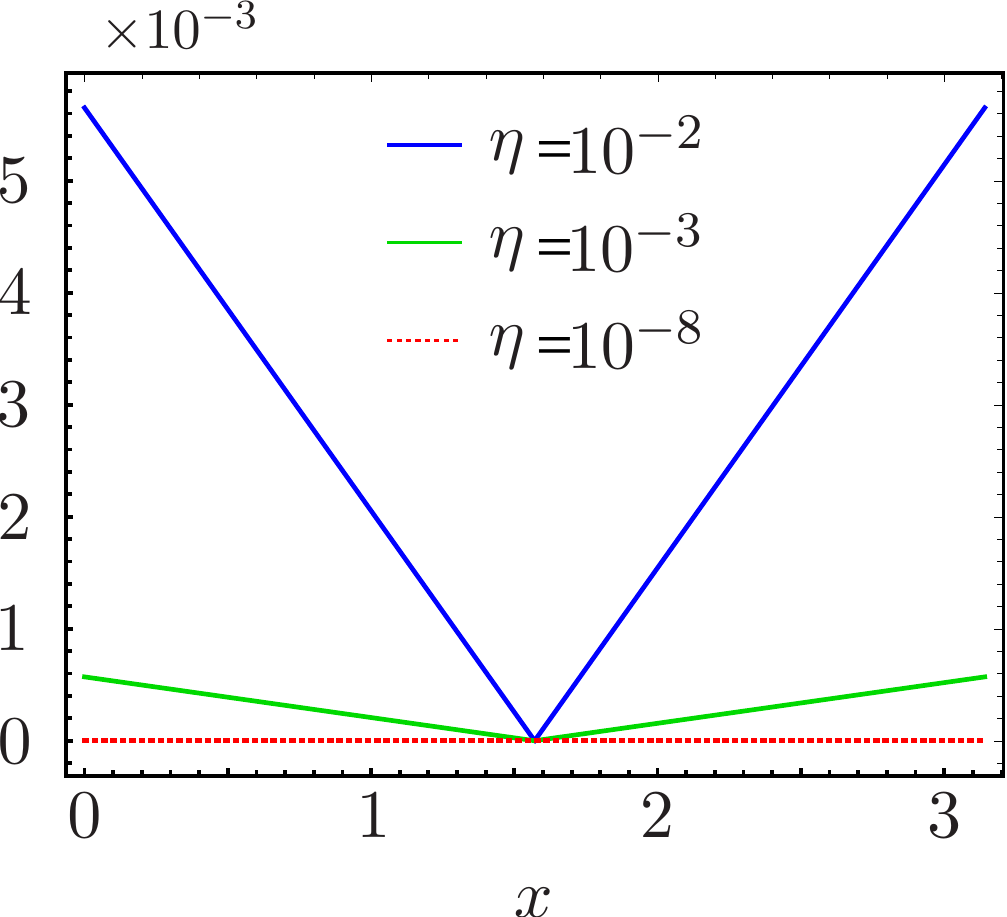} 
\end{center}
\caption{The difference $|v(x)-w(x)|$ vs. $x$ in the fluid domain $\Omega_f$ for $m=1$ and $\alpha=1$ for different values of the penalization parameter $\eta =10^{-2},\ 10^{-3}$ and $10^{-8}$. 
Here $w$ and $v$ are given by Eqs. (\ref{eq2-3}) and (\ref{eq2-9}), respectively.
} \label{figError}
\end{figure}

Now we solve Eq. (\ref{eq2-4}) analytically, and obtain
\begin{equation}
v(x)= \left \{
\begin{array}{ll}
\cos{mx}+A_1 x + A_2 & {\rm for} \quad 0 \leq x \leq \pi , \\
B_1x+B_2            & {\rm for} \quad \pi \leq x < 2\pi.
\end{array}
\right. \label{eq2-9}
\end{equation}
Using Eq. (\ref{eq2-8}) and the $C^0$-continuity of the solutions at the interface of $\Omega_f$ and $\Omega_s$, i.e.,
$v(0^{+})=v(2 \pi^{-})$ and $v(\pi^{+})=v(\pi^{-})$,
we can determine three out of the four constants.
\begin{eqnarray}
& & A_1 = \frac{\alpha}{1+\eta}+\eta \frac{1-(-1)^m}{\pi(1+\eta)}, \\
& & A_2 = 2\pi \left\{ -\frac{\alpha}{1+\eta}+\frac{1-(-1)^m}{\pi(1+\eta)} \right\} -1+B_2, \\
& & B_1 =  -\frac{\alpha}{1+\eta} + \frac{1-(-1)^m}{\pi(1+\eta)}.
\end{eqnarray}
Under the condition $\int_{\Omega_f} v (x) d x=0$, we get
\begin{eqnarray}
B_2 = \frac{3\pi}{2}\frac{\alpha}{1+\eta} -\frac{1-(-1)^m}{1+\eta}\left( \frac{\eta}{2} + 2 \right)+1 . 
\end{eqnarray}
In the limit of $\eta \rightarrow 0$, 
the penalized solution $v(x)$, Eq. (\ref{eq2-9}), converges towards the exact solution $w(x)$, Eq. (\ref{eq2-3}), of the non-penalized problem. 
For even $m$ and $\alpha=0$, $v(x)$ is identical to $w(x)$.
Fig. \ref{fig2-1} shows $v(x)$ at $\eta=10^{-1}$ and $10^{-8}$ for $m=1$ together with $w(x)$.
It is seen that for decreasing value of $\eta$, $v(x)$ perfectly matches $w(x)$.
For the first derivative $d_x v$, we observe that no penalization boundary layer is present, 
which is in good agreement with the findings in Ref. \cite{KoNgSc}. 
The second derivative $d_x^2 v$ confirms the absence of the penalization boundary layer.

The penalization error $\left| v(x)-w(x) \right|$ ($= \left| (A_1-\alpha)x+A_2-A_0 \right| $ in $\Omega_f$) is $O(\eta)$,
because $A_1=\alpha +O(\eta)$, $A_2=A_0+O(\eta)$ and $A_0=-\pi \alpha/2$.
We also find that $\left| d_x v -d_x w \right|$ is constant and we have again an error of $O(\eta)$.
Fig. \ref{figError} shows that the difference $\left| v(x) -w(x) \right|$ is symmetric with respect to $x=\pi/2$ and takes maximum values at the boundaries of the fluid domain, i.e., at $x=0$ or $x=\pi$.
Figs. \ref{fig2-1} and \ref{figError} suggest that
for the largest value of $\eta$, i.e., $\eta=10^{-1}$,
the penalization approach does not well approximate the boundary condition of the non-penalized problem.

\begin{figure}[tb]
\begin{center}
\includegraphics[width=7cm,keepaspectratio]{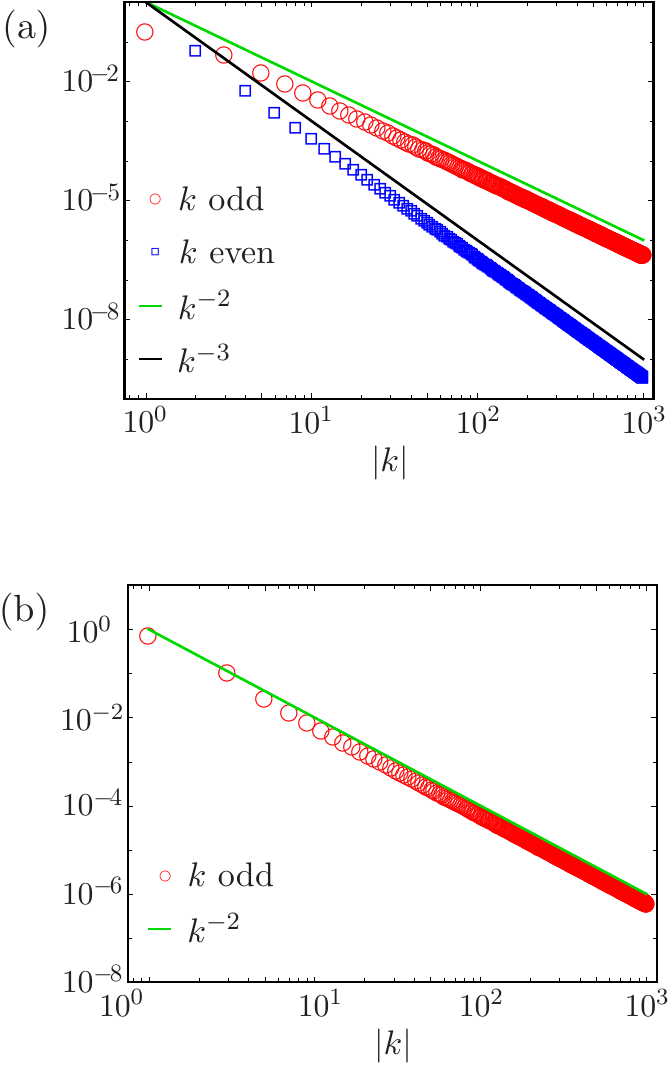}
\end{center}
\caption{Decay of the modulus of the Fourier coefficients $|{\hat v}_k|$ at $\eta = 10^{-2}$ for (a) $m=1$ and (b) $m=2$.
In Fig. \ref{fig2-2} (b), $|{\hat v}_k|$ is plotted only for odd $k$. }
\label{fig2-2}
\end{figure}

To study its regularity, we also analyze the Fourier coefficients ${\hat v}_k$ of the penalized solution $v(x)$, Eq. (\ref{eq2-9}),
where ${\hat v}_k$ is defined by 
\begin{equation}
{\hat v}_k = \frac{1}{2\pi} \int_0^{2\pi} v(x) \exp( -{\mathrm{i}} kx) dx.
\end{equation}
The Fourier coefficients ${\hat v}_k$ for odd $m$ can be computed explicitly and are given by 
\begin{eqnarray}
{\hat v}_k= \left \{
\begin{array}{ll}
\displaystyle{\frac{{\mathrm i}}{\pi}\frac{m^2}{k(m^2-k^2)}} & {\rm if} \ k \ {\rm even},  \\
& \\
\displaystyle{\frac{2}{\pi^2m^2} \frac{1-\eta}{1+\eta}-\frac{2\alpha}{\pi m^2}\frac{1}{1+\eta} +\frac{1}{4}}  & {\rm if} \ k \ {\rm odd \ and} \ k=\pm m ,\\
& \\
\displaystyle{\frac{2}{\pi^2k^2} \frac{1-\eta}{1+\eta}-\frac{2\alpha}{\pi k^2}\frac{1}{1+\eta}}  & {\rm if} \ k \ {\rm odd \ and} \ k \neq \pm m ,
\end{array}
\right. 
\label{eq2-11}
\end{eqnarray}
while the Fourier coefficients for even $m$ are given by
\begin{eqnarray}
{\hat v}_k= \left \{
\begin{array}{ll}
\displaystyle{-\frac{2\alpha}{\pi k^2}\frac{1}{1+\eta}+\frac{{\mathrm i}}{\pi}\frac{m^2}{k(m^2-k^2)}} & {\rm if} \ k \ {\rm odd},  \\
& \\
\displaystyle{\frac{1}{4}} & {\rm if} \ k \ {\rm even \ and} \ k=\pm m ,\\
& \\
0  & {\rm if} \ k \ {\rm even \ and} \ k \neq \pm m .
\end{array}
\right.
\label{eq2-12}
\end{eqnarray}
For sufficiently large $|k|$ in the sense that $|k| \gg |m|$, Eqs. (\ref{eq2-11}) and (\ref{eq2-12}) show
that $|{\hat v}_k|$ decays proportional to $k^{-2}$ for odd wave numbers $k$,
whereas it decays as $O(k^{-3})$ only for even $k$ and odd $m$.
Fig. \ref{fig2-2} shows the decay of $|{\hat v}_k|$ for $\eta=10^{-2}$.
There is no intermediate region at low $k$, 
as it is the case for the VP representation of homogeneous Neumann boundary conditions \cite{KoNgSc}. 
The absence of the region is in contrast to what is observed for the VP representation of the Dirichlet boundary condition \cite{Romain},
and is thus attributed to the absence of any boundary layer due to the penalization of Neumann boundary conditions.

\subsection{Discretization error of the second order finite-difference scheme} \label{sec2-2}
We discretize the penalized equation, 
using the second order central finite-differences of the first derivative of $\phi(x)$
given by 
\begin{equation}
d_x \phi (x) = \frac{\phi(x+ h/2)-\phi(x- h/2) }{h} +O(h^2) , \label{2ndfd}
\end{equation}
and second order interpolation of $\phi$,  
\begin{equation}
\phi(x) = \frac{\phi(x + h/2)+\phi(x - h/2)}{2} +O(h^2) , \label{2ndint}
\end{equation}
where $\phi=v, d_x v, \, \theta d_x v $, $x=x_i, x_i\pm h/2$, $x_i=i h$ $(i=0,\cdots,N-1)$, $h \, (=2\pi/N)$ is the grid width, and $N$ is the number of grid points.

\begin{figure}[tb]
\begin{center}
\includegraphics[width=7cm,keepaspectratio]{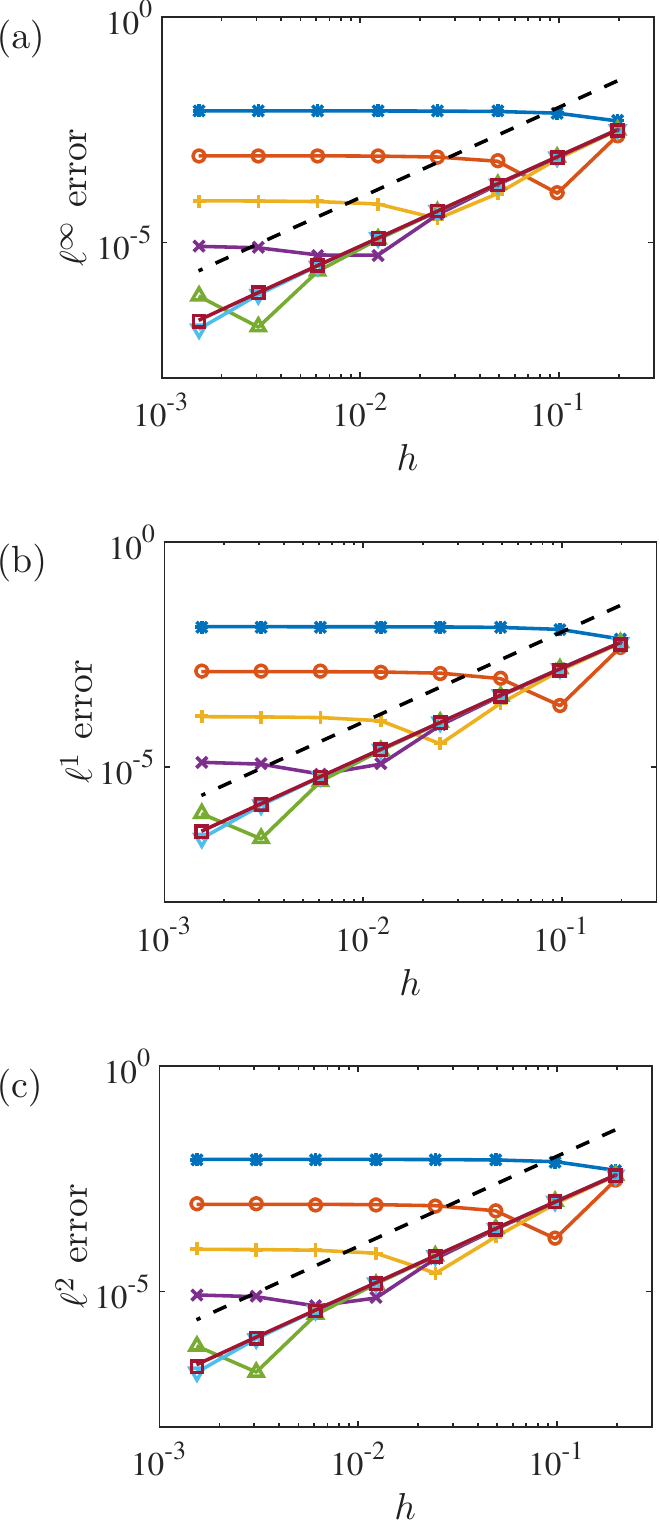}
\end{center}
\caption{The errors between the numerical solution of Eq. (\ref{eq2-4}) and the exact solution (\ref{eq2-3}) vs. $h$ for $\alpha = 0.1$ at $\eta =10^{-2},10^{-3},10^{-4},10^{-5},10^{-6},10^{-7}$ and $10^{-8}$, 
which are respectively denoted by the line with $*$, $\circ$, $+$, $\times$, $\triangle$, $\triangledown$ and $\square$.
The errors are measured by (a) the $\ell^{\infty}$-norm, (b) the $\ell^{1}$-norm, and (c) the $\ell^{2}$-norm.
The dashed lines show the $O(h^2)$ decay.}
\label{fig2-4-1}
\end{figure}
\begin{figure}[tb]
\begin{center}
\includegraphics[width=7cm,keepaspectratio]{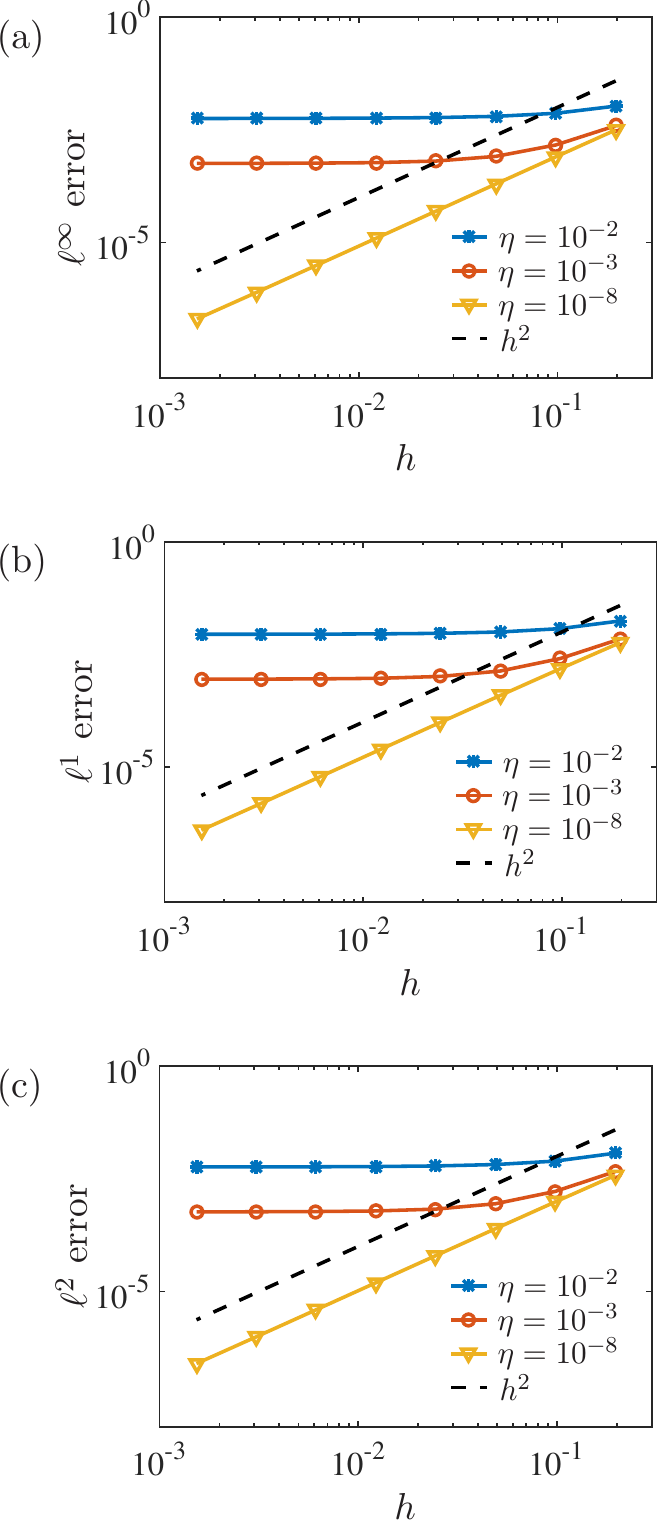}
\end{center}
\caption{The errors between the numerical solution of Eq. (\ref{eq2-4}) and the exact solution (\ref{eq2-3}) vs. $h$ for $\alpha = 1$ at $\eta =10^{-2},10^{-3}$ and $10^{-8}$.
The errors are measured by (a) the $\ell^{\infty}$-norm, (b) the $\ell^{1}$-norm, and (c) the $\ell^{2}$-norm.}
\label{fig2-4-2}
\end{figure}

We then analyze the discretization and penalization errors.
In Figs. \ref{fig2-4-1} and \ref{fig2-4-2}, 
we plot the errors between the numerical solution of the penalized equation (\ref{eq2-4}) and the exact solution (\ref{eq2-3}) of the non-penalized problem, Eqs. (\ref{eq2-0})--(\ref{eq2-2}), as a function of $h$. 
The errors are measured in $\ell^\infty$-, $\ell^1$- and $\ell^2$-norms, and are computed only in the fluid domain $\Omega_f$. 
The grid points $x_0$ and $x_{N/2}$ are on the boundary of $\Omega_f$.

For small $\alpha \, (=0.1)$, Fig. \ref{fig2-4-1} shows that the error measured in each norm has a pronounced minimum, 
which corresponds to the optimal value of $h$, for a given penalization parameter $\eta$. 
By decreasing $h$ below this value, 
we see that the errors increase, and eventually saturate. 
This saturation is due to the penalization error.
For smaller $\eta$, 
the errors decay approximately as $O(h^2)$ with decreasing $h$. 
The minimum value and the optimal value of $h$ decrease, 
as $\eta$ becomes smaller.
These behaviors are the same as those obtained by Ref. \cite{KoNgSc} for the penalized problem of the homogeneous Neumann boundary condition, i.e., $\alpha=0$.
It can be seen that
for $h \approx 10^{-2}$, the error has its minimum for $\eta=10^{-5}$, for the three norms shown in Fig. \ref{fig2-4-1}.

However, there is no pronounced minimum for larger values $\alpha \, (=1)$.
Fig. \ref{fig2-4-2} shows that, for larger $\eta \, (=10^{-2}) $,
the errors for all considered norms are almost independent of $h$.
For $\eta=10^{-3}$, the errors decay approximately as $O(h^2)$ with decreasing $h$ and then saturate.
The level of the saturation becomes smaller, as $\eta$ decreases.
The $O(h^2)$ convergence is again observed at sufficiently small $\eta$.
For much larger $\alpha \,(=10)$, 
we observe the behavior similar to the case of $\alpha=1$ (figure omitted).

It is suggested in Section \ref{sec2-3} and Appendix \ref{appeA} that this $O(h^2)$ convergence is attributed to the cancellation of the errors of $O(h)$.
The cancellation is due to the symmetry of Eq. (\ref{eq2-4}) with respect to $x=\pi/2$.

The discrete eigenvalue problem of the penalized Laplace operator results in the same as the problem in the homogeneous case,
computed in Ref. \cite{KoNgSc}, 
that shows two different behaviors in the spectrum of the penalized operator.
The eigenvalues of one part of the spectrum converge to the eigenvalues of the non-penalized Laplace operator.
For the other part, 
corresponding to the small eigenvalues, 
these values depend on $\eta$ and vanish at $\eta \rightarrow 0$.
Let $F$ be the discretized and penalized Laplace operator which is obtained by application of the second order finite-difference scheme, Eqs. (\ref{2ndfd}) and (\ref{2ndint}), to $- d_x ( \theta d_x)$.
The corresponding eigenvalue problem has the form $F {\bm v} = {\bm b}$.
The circulant matrix $F$ can be diagonalized applying Fourier series expansion, 
and we have ${\mathrm{ker}} F = 1$. 
Hence it is not invertible. 
To remove this kernel, we use the constraint, a discretized version of $\int_{\Omega_f} v (x) dx =0$ using the midpoint rule.
Under periodic boundary conditions,
the eigenfunctions can be represented by the Fourier series.
Thus, all of the non-zero eigenvalues are positive.

\subsection{1D Poisson equation with different inhomogeneous Neumann boundary conditions}  \label{sec2-3}
In this section, we study a VP representation for two different values of inhomogeneous Neumann boundary conditions.
We consider again a 1D Poisson equation Eq. (\ref{eq2-0}), $-d_x^2 w(x) = f(x)$, with the source term given by
\begin{eqnarray}
f(x) = m^2 \sin(m x)  \label{eq2-14-2}
\end{eqnarray}
for odd $m$ in the domain $\Omega_f=\{ x \ | \ 0 <x <\pi\}$ and imposing
\begin{eqnarray}
d_x w(0) = \alpha+m, \quad {\mathrm{and}} \quad d_x w(\pi)= \alpha - m . \label{eq2-15}
\end{eqnarray}
The compatibility condition is fulfilled,
since we obtain $ -[d_x w]_0^\pi = 2m $ and $\int_0^\pi f (x) dx =2m $.
Imposing that $\int_{\Omega_f} w (x) dx =0$, 
we find that the exact solution is given by
\begin{eqnarray}
w(x) = \sin(mx) + \alpha x -\frac{2}{m \pi} - \frac{\pi}{2}\alpha.  \label{eq2-16}
\end{eqnarray}

The VP representation of Eqs. (\ref{eq2-0}) and (\ref{eq2-15}) reads 
\begin{eqnarray}
-  d_x ( \theta d_x v + \beta \chi ) = (1-\chi) f - \chi d_x \beta,
\label{eq2-17}
\end{eqnarray}
in the extended domain $\Omega=\{ x \ | \ 0 \leq x < 2\pi \}$ under $2\pi$-periodic boundary conditions,
where the mask function $\chi(x)$ is given by Eq. (\ref{eq2-7}).
Here we set $\beta(x)=\alpha + m\cos(m x)$, which satisfies $\beta(0)=\alpha + m$ and $\beta(\pi)=\alpha -m$.
The second term on the right-hand side of Eq. (\ref{eq2-17}), $- \chi d_x \beta$, is added 
such that Eq. (\ref{eq2-17}) reduces to $\eta d_x^2 v=0$ in the solid domain $\Omega_s \, (=\Omega \backslash {\bar \Omega}_f)$,
and thus the penalized solution converges towards the non-penalized one (\ref{eq2-16}) for $\eta \rightarrow 0$.
The penalized representation (\ref{eq2-17}) is again based on the continuity of the flux though the non-penalized boundary at $x=0 \, (=2\pi),\pi$.
The flux at the interface satisfies
\begin{eqnarray}
d_x v = \eta d_x v + \beta \quad {\rm at} \quad x=0 \, (=2\pi), \ \pi . \label{eq2-19}
\end{eqnarray}

\begin{figure}[tb]
\begin{center}
\includegraphics[width=7cm,keepaspectratio]{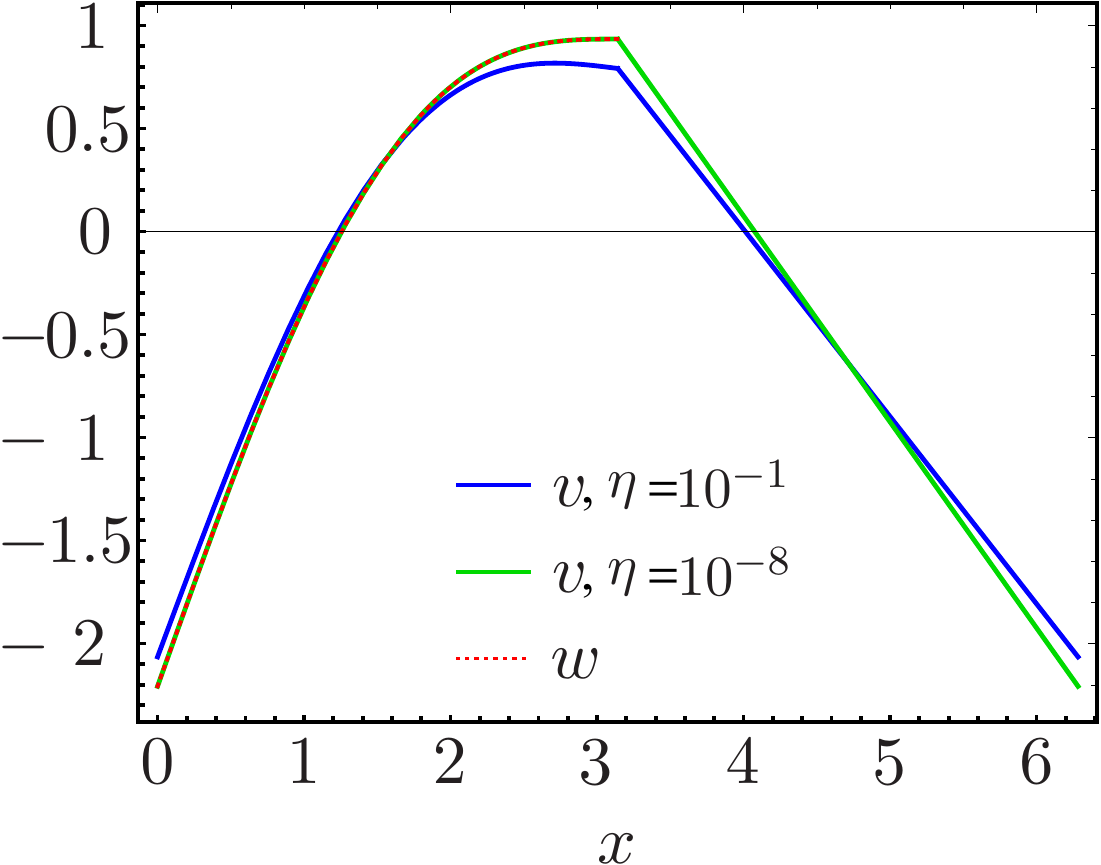}
\end{center}
\caption{The exact solution $w$, Eq. (\ref{eq2-16}), in the fluid domain $\Omega_f$ at $\alpha =1$ and $m=1$ (red line). 
The analytical solution $v$, Eq. (\ref{eq2-20}), in the extended domain $\Omega$ at $\eta = 10^{-1}$ (blue line) and $\eta = 10^{-8}$ (green line).}
\label{fig2-x}
\end{figure}

\begin{figure}[tb]
\begin{center}
\includegraphics[width=7cm,keepaspectratio]{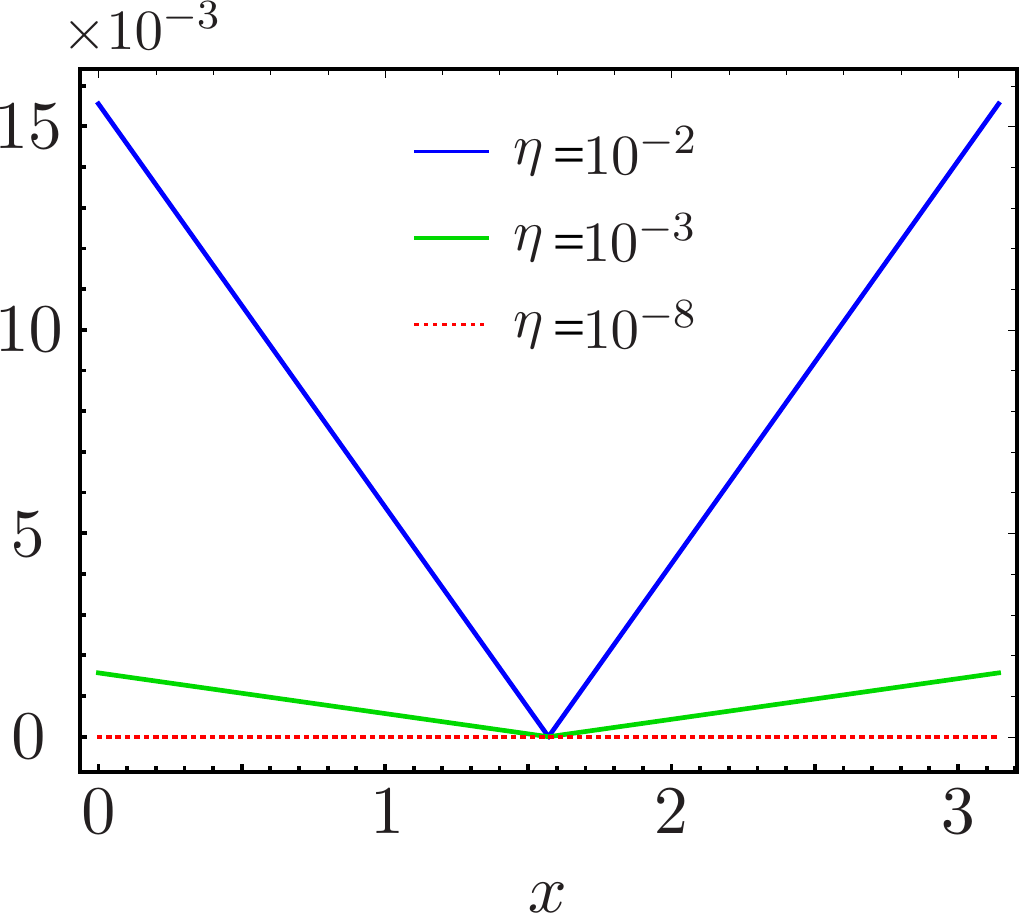}
\end{center}
\caption{The difference $|v(x)-w(x)|$ vs. $x$ in $\Omega_f$ for $m=1$ and $\alpha = 1$ at $\eta =10^{-2},\ 10^{-3}$ and $10^{-8}$,
where $w$ and $v$ are given by Eqs. (\ref{eq2-16}) and (\ref{eq2-20}), respectively.}
\label{fig2-z}
\end{figure}

We obtain the penalized solution of Eq. (\ref{eq2-17}); 
\begin{eqnarray}
v(x) =\left \{
\begin{array}{ll}
\sin(m x)+A_1' x +A_2' & {\rm for} \quad 0 \leq x \leq \pi , \\
B_1' x+B_2'          & {\rm for} \quad \pi \leq x < 2\pi .
\end{array} \right. \label{eq2-20}
\end{eqnarray}
The coefficients can be determined in the same way for Eq. (\ref{eq2-9}), 
and we obtain
\begin{eqnarray}
A_1' \!\!\!&=&\!\!\! \frac{\alpha}{1+\eta}, \quad  A_2'= - \frac{\pi}{2} \frac{\alpha }{1+\eta} - \frac{2}{m\pi} ,
\label{eq2-21-0} \\
B_1' \!\!\!&=&\!\!\! -\frac{\alpha}{1+\eta}, \quad  B_2'= \frac{3 \pi}{2} \frac{\alpha}{1+\eta} - \frac{2}{m\pi}.
\label{eq2-21}
\end{eqnarray}
The penalized solution $v$, Eq. (\ref{eq2-20}), results in the non-penalized solution $w$, Eq. (\ref{eq2-16}), for $\eta \rightarrow 0$.
We find $|v-w|=O(\eta)$ in the fluid domain $\Omega_f$, using Eq. (\ref{eq2-21-0}). 
Again it is found that $\left| d_x v - d_x w \right|$ is constant and of $O(\eta)$. 
Equation (\ref{eq2-20}) is independent of $\eta$ for $\alpha=0$.
Fig. \ref{fig2-x} presents the comparison of the penalized solution $v$ with the non-penalized solution $w$ for $\alpha=1$ and $m=1$ using different values of $\eta$.
It is seen that they excellently agree with each other for $\eta=10^{-8}$.
The difference $|v-w|$, which characterizes the penalization error at each position $x$, is shown in Fig. \ref{fig2-z}.
For each $\eta$, it is symmetric with respect to $x=\pi/2$.
The maximum value is located at the interface $x=0$ or $x=\pi$,
and decreases as $\eta$ becomes smaller.

\begin{figure}[tb]
\begin{center}
\includegraphics[width=7cm,keepaspectratio]{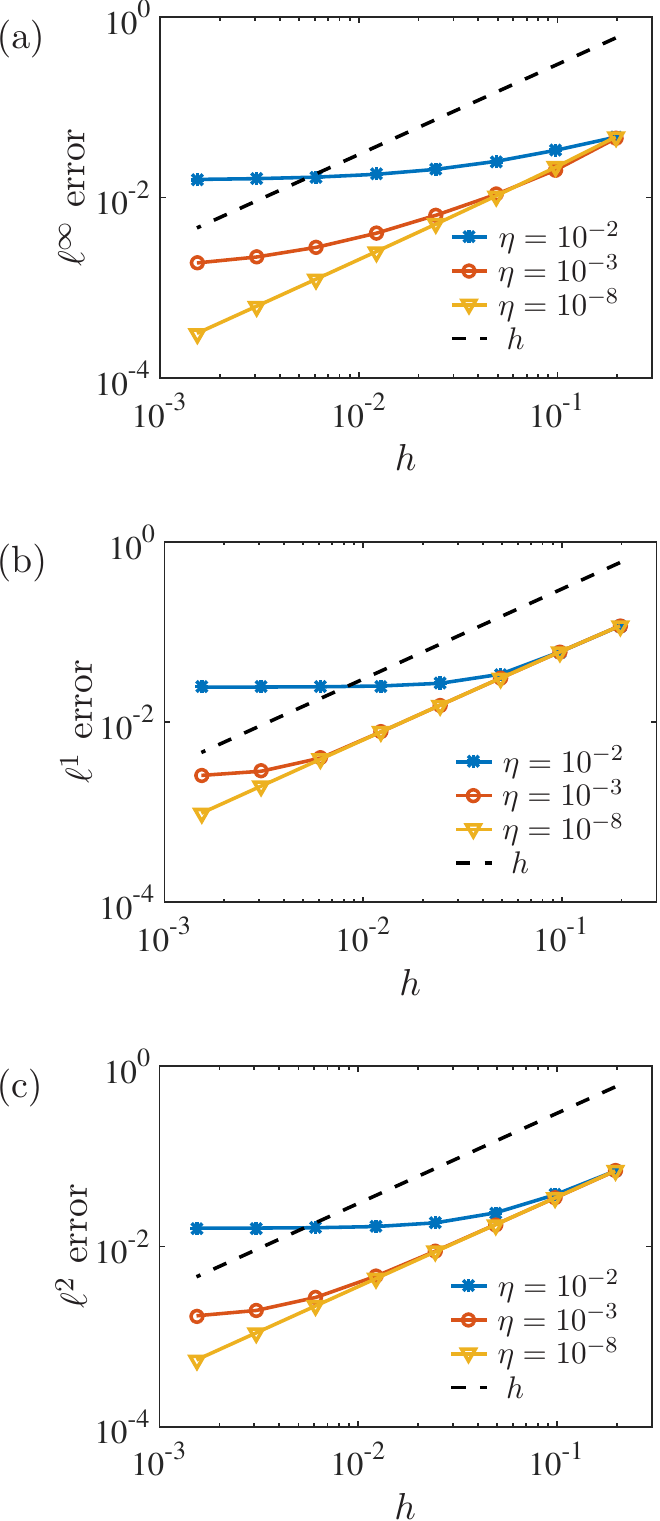}
\end{center}
\caption{The errors between the numerical solution of Eq. (\ref{eq2-19}) and the exact solution (\ref{eq2-16}) vs. $h$ for $m=1$ and $\alpha = 1$ at $\eta =10^{-2},10^{-3}$ and $10^{-8}$.
The errors are computed in $\Omega_f$, 
and measured by (a) the $\ell^{\infty}$-norm, (b) the $\ell^{1}$-norm and (c) the $\ell^{2}$-norm.
The dashed lines show the $O(h)$ decay.
}
\label{fig2-y}
\end{figure}

The numerical solutions of the non-penalized problem, Eqs. (\ref{eq2-0}), (\ref{eq2-14-2}) and (\ref{eq2-15}), and the penalized equation (\ref{eq2-17}) are obtained by using the second order finite-differences and interpolation in Eqs. (\ref{2ndfd}) and (\ref{2ndint}).
The errors measured by the different norms are plotted in Fig. \ref{fig2-y}.
The errors decay with decreasing $\eta$ and the grid width $h (=2\pi/N)$.
For smaller 
$\eta \,(= 10^{-3}$, $10^{-8})$, 
we observe that the errors exhibit first order convergence in terms of $h$.
This is in contrast to the case of Section \ref{sec2-1} where identical values of the inhomogeneous boundary conditions are considered on the left and the right boundary.
This $O(h)$ convergence is thus due to the different values of $\beta(x)$ at the interfaces $x=0,\pi$, i.e., $\beta(0) \ne \beta(\pi)$ (see Appendix \ref{appeA}).

The developed VP method for inhomogeneous Neumann boundary conditions can be also applied to time-dependent problems.
A convergence issue in terms of $\eta$ and $h$ is discussed in Appendix \ref{appB} for a 1D penalized thermal diffusion equation.

%
\section{Application to a 2D Poisson equation} \label{sec3}
Now, we extend the VP representation for inhomogeneous Neumann boundary conditions to 2D domains.
Let us consider a Poisson equation in the 2D fluid domain $\Omega_f=\{(x,y) \ | \ \pi/2 < x < 3\pi/2 \ {\rm and} \ \pi/2 < y < 3\pi/2 \}$;
\begin{equation}
-\nabla^2 w({\bm x}) = f({\bm x}), \label{eq3-1}
\end{equation}
imposing inhomogeneous Neumann boundary conditions ${\bm n}{\bm \cdot} \nabla w   = {\bm n}{\bm \cdot} {\bm \alpha} $ at the boundary of $\Omega$, i.e., $x=\pi/2, 3\pi/2$ for $\pi/2 \le y \le 3\pi/2$, and $y=\pi/2, 3\pi/2$ for $\pi/2 \le x \le 3\pi/2$,
where ${\bm x}=(x,y)$, $\nabla=(\partial_x,\partial_y )=(\partial/\partial x, \partial/\partial y )$, 
${\bm \alpha}=(\alpha_x, \alpha_y)$, $\alpha_x$ and $\alpha_y$ are constant,
and ${\bm n}$ is the outward pointing unit normal vector of $\Omega_f$.
We set the source function as $f({\bm x})=5\sin{x}\cos{2y}$, 
and verified the compatibility condition of Eq. (\ref{eq3-1}).
Then we obtain an exact solution;
\begin{eqnarray}
w({\bm x}) = \sin{x} \cos{2y} + \alpha_x x + \alpha_y y -(\alpha_x + \alpha_y)\pi. \label{eq3-3} 
\end{eqnarray}
We imposed the condition that $\int \!\! \int_{\Omega_f} w({\bm x}) d {\bm x} = 0$, to determine the solution uniquely.

A VP representation of Eq. (\ref{eq3-1}) reads
\begin{equation}
-\nabla \cdot ( \theta \nabla v + \chi {\bm \alpha} ) = (1-\chi) f , \label{eq3-4}
\end{equation}
which may be solved in the double periodic domain $\Omega=\{ (x,y) \ | \ 0 \le x < 2\pi,  0 \le y < 2\pi \}$.
The mask function $\chi({\bm x})$ is defined by
\begin{equation}
\chi({\bm x}) = \left\{
\begin{array}{ll}
0 & {\rm if} \,\,\,  {\bm x} \in \Omega_f ,\\
1/4 & {\rm if} \,\,\, {\bm x} = ( \pi/2, \pi/2 ), (\pi/2, 3\pi/2 ), \\
    & \quad \quad  \,\,\,\,\,\,\, ( 3\pi/2, \pi/2 ), ( 3\pi/2, 3\pi/2 ),  \\
1/2 & {\rm if} \,\,\, x=\pi/2,\, 3\pi/2 \,\,{\mathrm{for}}\,\, \pi/2<y<3\pi/2, \,\, {\mathrm{or}} \,\, y=\pi/2, \, 3\pi/2 \,\,{\mathrm{for}}\,\, \pi/2<x<3\pi/2 \\
1 & {\rm otherwise}  .
\end{array}
\right.
\end{equation}

\begin{figure}[tb]
\begin{center}
\includegraphics[width=7cm,keepaspectratio]{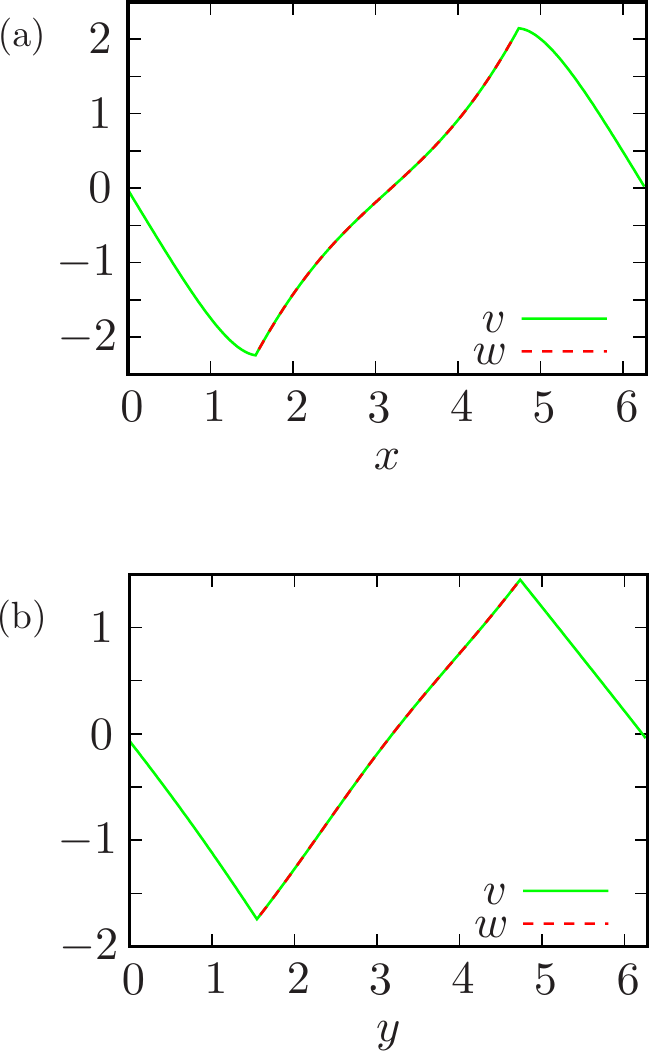}
\end{center}
\caption{The numerical solution $v$ of Eq. (\ref{eq3-4}) is compared with the exact solution $w$, Eq. (\ref{eq3-3}), at ${\bm \alpha}=(2,1)$ and $\eta=10^{-8}$:
(a) $w$ and $v$ vs. $x$ at $y=3.09$, and (b) $w$ and $v$ vs. $y$ at $x=3.09$.}
\label{fig3-1}
\end{figure}

Imposing $\int \!\! \int_{\Omega_f} v({\bm x}) d {\bm x}=0$, 
we select the penalized solution of Eq. (\ref{eq3-4}) that converges to the exact non-penalized solution (\ref{eq3-3}) for $\eta \rightarrow 0$.
Fig. \ref{fig3-1} illustrates the numerical solution of Eq. (\ref{eq3-4}) at $\eta=10^{-8}$ for ${\bm \alpha}=(2,1)$, 
which is obtained by the second order finite-differences and interpolation in Eqs. (\ref{2ndfd}) and (\ref{2ndint}), 
using $256$ grid points in each Cartesian direction.
We observe that the penalized solution $v({\bm x})$ is in excellent agreement with the exact solution (\ref{eq3-3}) in $\Omega_f$.

\begin{figure}[tb]
\begin{center}
\includegraphics[width=7cm,keepaspectratio]{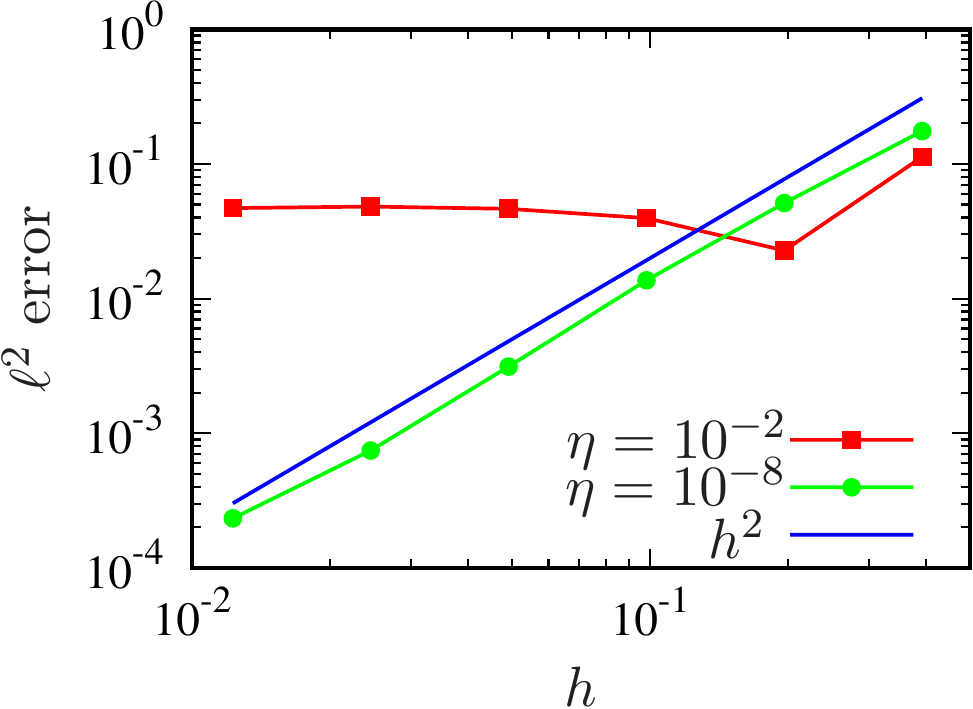}
\end{center}
\caption{The $h$-dependence of the $\ell^2$ errors between the numerical solution of Eq. (\ref{eq3-4}) and the exact solution (\ref{eq3-3}) at $\eta=10^{-2}$ and $10^{-8}$. 
The errors are computed in the fluid domain $\Omega_f$.
The blue line shows the $O(h^2)$ decay.}\label{fig3-2}
\end{figure}

Fig. \ref{fig3-2} shows the dependence of the $\ell^2$ error between the numerical solution $v$ and the exact solution $w$, Eq. (\ref{eq3-3}), in $\Omega_f$ on the grid width $h$ for $\eta=10^{-2}$ and $10^{-8}$. 
For $\eta=10^{-2}$, the $\ell^2$ error has the minimum in terms of $h$.
The error becomes almost constant for smaller values of $h \, (\lesssim 0.1)$.
In contrast, for much smaller $\eta \, (=10^{-8})$,
it is seen that the error monotonically decays with decreasing $h$, and behaves approximately as $O(h^2)$.
This $O(h^2)$ convergence is consistent with what was found in Fig. \ref{fig2-4-2} for the 1D penalized Poisson equation,
as well as for the 2D penalized Poisson equation with homogeneous Neumann boundary conditions in a square domain \cite{KoNgSc}.
The errors measured by the $\ell^1$- and $\ell^\infty$-norms are omitted, 
because the errors show the same behaviors as those of the $\ell^2$ error.


Next, let us consider a 2D Poisson equation (\ref{eq3-1}) with inhomogeneous Neumann boundary conditions in an annular fluid domain $\Omega_f=\{(x,y) \ | \ \pi/4 < r < 3\pi/4 \} $, 
where $r=\sqrt{(x-\pi)^2+(y-\pi)^2}$.
We set the source term $f({\bm x})$ as  
\begin{equation}
f({\bm x}) = f(r) = 16\cos{4r}+\frac{4\sin{4r}}{r} ,\label{Y0}
\end{equation}
and impose inhomogeneous Neumann boundary conditions given by
\begin{equation}
\frac{\partial w}{\partial r} = 3\alpha \,\,\, {\mathrm{at}} \,\,\, r=\frac{\pi}{4}, 
\,\,\, {\mathrm{and}} \,\,\, \frac{\partial w}{\partial r} = \alpha \,\,\, {\mathrm{at}} \,\,\, r= \frac{3}{4} \pi.  \label{Y1}
\end{equation}
The compatibility condition $\int_{\pi/4}^{3\pi/4} r f(r) dr=0$ is fulfilled.
Noting $\nabla^2 w =(1/r) d_r (r d_r)$ and $w({\bm x}) =w(r)$, 
we obtain the exact solution;
\begin{eqnarray}
w(r)=\cos{4r}+\frac{3}{4}\alpha \pi \log{r}+C , \label{exactC}
\end{eqnarray}
where 
\begin{eqnarray}
C= - \frac{3}{32}\alpha\pi\left( 9\log{\frac{3}{4}\pi} - \log{\frac{\pi}{4}} -4 \right),
\end{eqnarray}
imposing the condition that $\int_{\pi/4}^{3\pi/4} r w(r) d r  =0$.

The penalized equation of this problem reads
\begin{equation}
-\nabla \cdot ( \theta \nabla v + \chi {\bm \beta} ) = (1-\chi) f  - \chi \nabla \cdot {\bm \beta}, \label{eqCP}
\end{equation}
where $\Omega_f$ is imbedded into the $2\pi$ double periodic domain $\Omega=\{(x,y)| 0\le x < 2\pi, \, 0 \le y < 2\pi \}$,
${\bm \beta}({\bm x}) =g({\bm x}) {\bm e}_r$, ${\bm e}_r=(x-\pi,y-\pi)/r$ and
\begin{eqnarray}
g({\bm x}) &=&\left\{
\begin{array}{ll}
\displaystyle{\alpha \left( \frac{4r}{3\pi} \right)^2 \left\{ 4\left(1 -\frac{r}{\pi} \right) \right\}^3 } &{\rm if} \quad 0\leq r \leq \pi ,\\
0 & {\rm otherwise} .
\end{array}
 \right. 
\end{eqnarray}
Note that ${\bm \beta}({\bm x})=3\alpha {\bm e}_r$ at $r =\pi/4$, and ${\bm \beta}({\bm x})=\alpha {\bm e}_r$ at $r=3\pi/4$,
which is consistent with Eq. (\ref{Y1}).
The mask function $\chi({\bm x})$ is given by
\begin{eqnarray}
\chi(x)= \left \{
\begin{array}{ll}
0   & {\rm if} \quad \pi/4 < r < 3\pi/4 , \\
1/2 & {\rm if} \quad r=\pi/4, \  3\pi/4 , \\
1   & {\rm otherwise}.
\end{array}
\right. 
\end{eqnarray}

\begin{figure}[tb]
\begin{center}
\includegraphics[width=7cm,keepaspectratio]{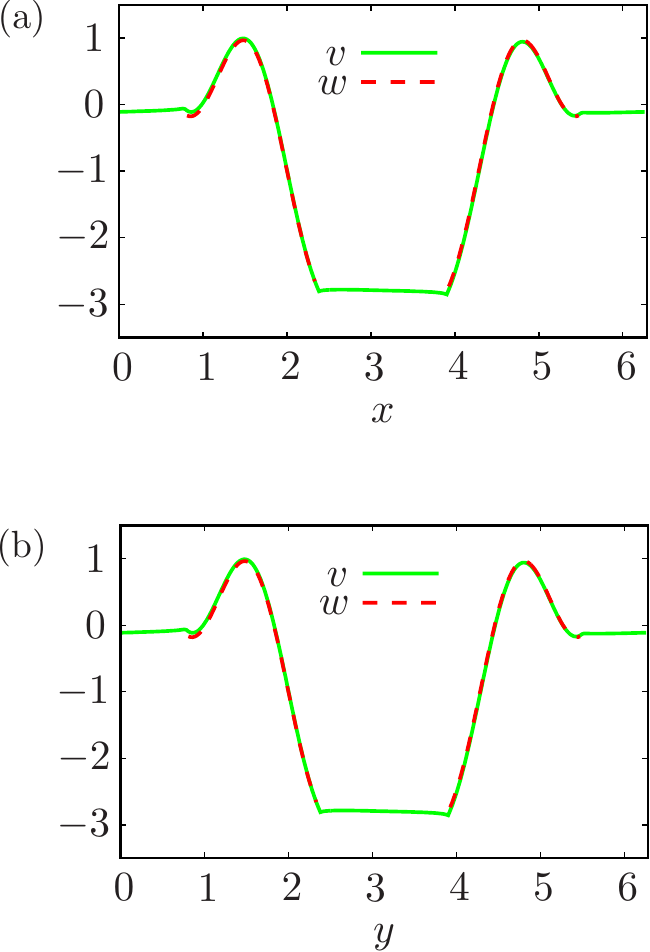}
\end{center}
\caption{1D cuts of the numerical solution $v$ of Eq. (\ref{eqCP}) and the exact solution $w$, Eq. (\ref{exactC}), at $\alpha=1$ and $\eta=10^{-8}$:
(a) $v$ and $w$ vs. $x$ at $y=3.12$, and (b) $v$ and $w$ vs. $y$ at $x=3.12$.}
\label{fig3-3}
\end{figure}

\begin{figure}[tb]
\begin{center}
\includegraphics[width=7cm,keepaspectratio]{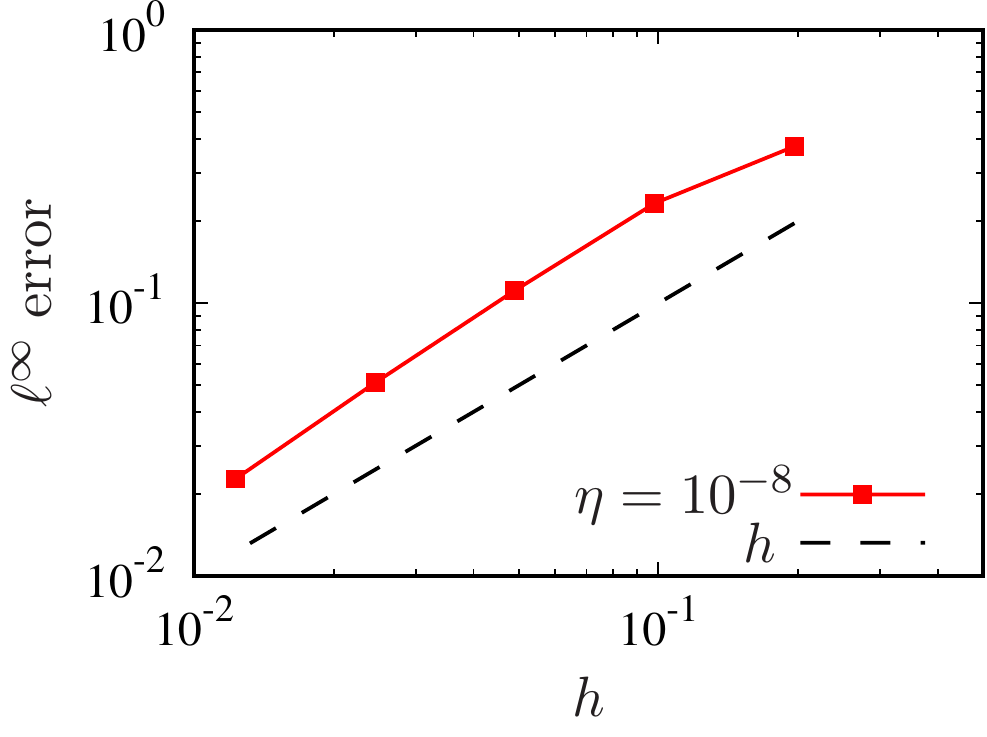}
\end{center}
\caption{The $h$-dependence of the $\ell^\infty$ error between the numerical solutions of Eq. (\ref{eq3-4}) and the exact solution (\ref{exactC}) at $\alpha=1$ and $\eta=10^{-8}$. The error is computed in $\Omega_f$.
The green line shows the $O(h)$ decay.}
\label{fig3-4}
\end{figure}

Fig. \ref{fig3-3} shows 1D cuts of the numerical solution $v$, 
which are compared with those of the exact solution $w$, Eq. (\ref{exactC}), of the non-penalized problem for $\alpha=1$ and $\eta=10^{-8}$.
The numerical solutions were computed in Cartesian coordinates using $256$ grid points in each direction.
It is seen that $v$ is in good agreement with $w$.
Fig. \ref{fig3-4} shows the dependence of the $\ell^\infty$ error between the numerical solution $v$ and the exact one $w$ as a function of the grid width $h$ in $\Omega_f$.
The error decays approximately as $O(h)$.
This $O(h)$ convergence is attributed to the following two errors.
One is the absence of grid points at the circular interfaces, i.e. for $r=\pi/4$ and $3\pi/4$ (see Ref. \cite{KoNgSc}).
The other is the use of the different boundary values for the inhomogeneous Neumann boundary conditions, as already discussed in Section \ref{sec2-3} and Appendix \ref{appeA}.


\section{Application to free convection in two dimensions} \label{sec5}
We now apply the VP method developed in this paper to an advection-diffusion equation with an inhomogeneous Neumann boundary condition and a Dirichlet boundary condition.  
The equation is coupled with the incompressible Navier--Stokes equations.
Numerical simulations were performed for 2D steady thermal flows in a concentric annulus that is heated from its inner cylinder. 
The VP method is examined by comparison with numerical results of Refs. \cite{Ren2013,Yoo2003}. 
\begin{figure}[tb]
\begin{center}
\includegraphics[width=7cm,keepaspectratio]{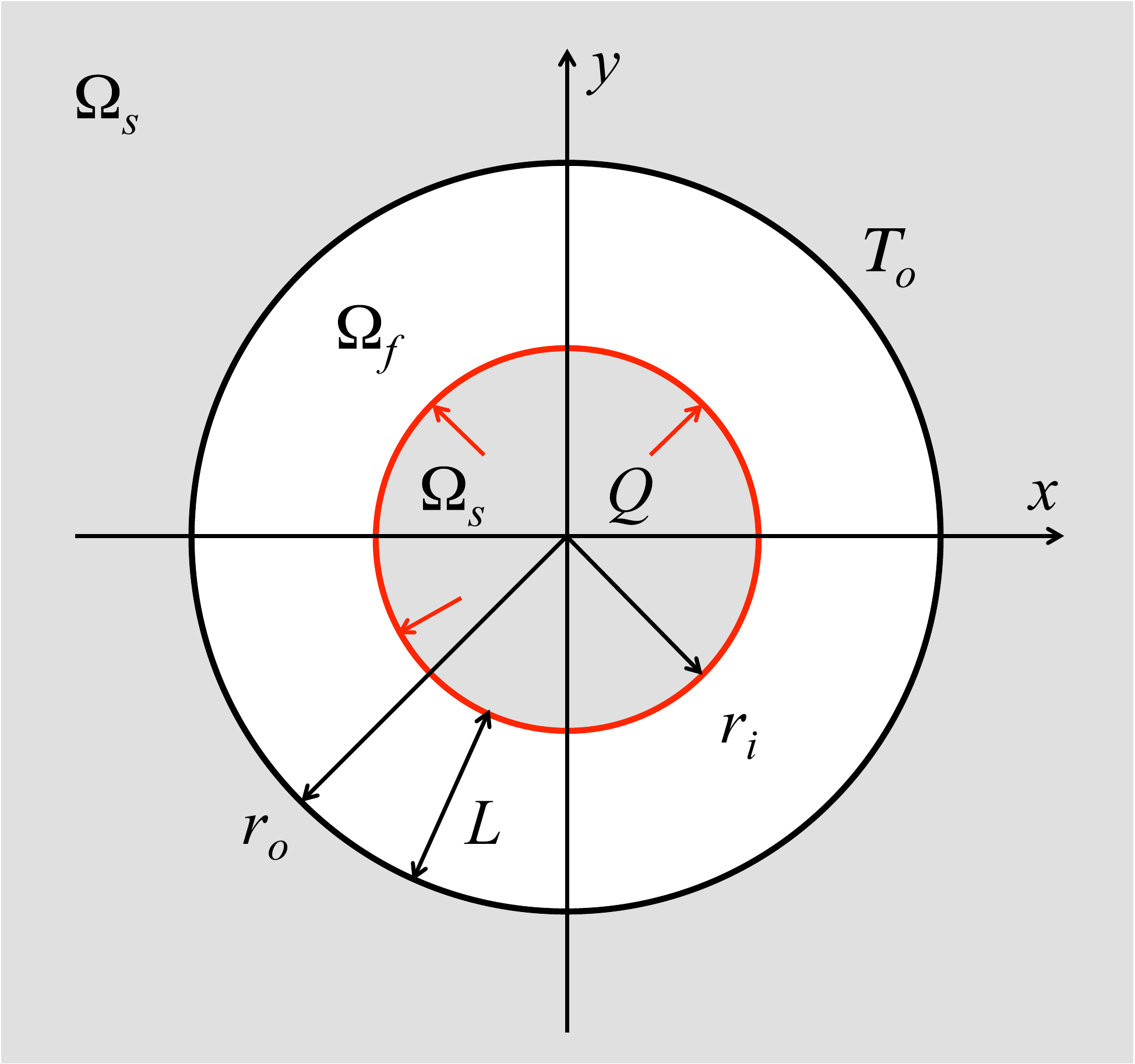}
\end{center}
\caption{Configuration for a concentric horizontal cylindrical annulus.
The gravity force is imposed as $(0,-g)$, where $g$ is the modulus of the gravity acceleration.}\label{fig5-0}
\end{figure}
%
\subsection{Volume penalized governing equations} \label{sec5-1}
We consider steady free convection of an incompressible fluid in a concentric annulus.
Its configuration is illustrated in Fig. \ref{fig5-0}, 
where $r_i$ and $r_o$ are respectively the inner and outer cylinder radius, 
$L=r_o-r_i$, $T_o$ is constant temperature of the outer cylinder, 
and $Q$ is a constant heat flux from the inner cylinder.
In this study, we set $r_i/L = 1$ and $r_o/L=2$, 
which was used in Refs. \cite{Ren2013,Yoo2003}. 
Let ${\bm u}$, $T$ and ${\bm x}$ be respectively the velocity, the temperature and the position.
The dimensionless velocity ${\bm u}^{*}$ is given by ${\bm u}^{*} = {\bm u} L/\kappa $, 
the dimensionless temperature $\Phi^{*}$ is denoted by $\Phi^{*} = (T - T_o) k/ (Q L)$, and ${\bm x}^{*} = {\bm x}/{L}$,
where $\kappa $ is the thermal diffusivity defined by $\kappa= k / (\rho C_p) $; $k$ is the thermal conductivity, $\rho$ is the density, and $C_p$ is the heat capacity at constant pressure.
We will drop the superscript $^*$ hereafter, unless otherwise stated.

The non-penalized governing equations under the Boussinesq approximation are given by
\begin{eqnarray}
& & \partial_t {\bm u} = -\nabla p -({\bm u}  {\bm \cdot} \nabla) {\bm u}  + {\rm Pr} \nabla^2 {\bm u} + {\rm Pr Ra} \Phi {\bm e}_y , \label{eq5-1-1} \\
& & \nabla {\bm \cdot} {\bm u} = 0 , \label{eq5-1-2}\\
& & \partial_t \Phi= -({\bm u} {\bm \cdot} \nabla) \Phi + \nabla^2 \Phi, \label{eq5-1-3}
\end{eqnarray}
where $p$ is the dimensionless pressure, $\partial_t = \partial /\partial t$, ${\bm e}_y=(0,1)$, 
${\mathrm{Pr}} = \nu / \kappa$ is the Prandtl number, 
${\mathrm{Ra}} = g \gamma QL^4/(k \kappa \nu ) $ is the Rayleigh number,
and $\gamma$ is the thermal expansion coefficient.
The fluid domain $\Omega_f$ is given by $\Omega_f=\{(x,y) \ | \ 1 < r < 2  \}$, 
where $r=\sqrt{x^2+y^2}$.
No-slip boundary conditions for velocity, ${\bm u}={\bm 0}$, are imposed at $r=1$ and $r=2$. 
The boundary conditions of $\Phi$ are given by
\begin{eqnarray}
{\bm n} {\bm \cdot} \nabla \Phi = 1 & \quad {\rm at} \quad r = 1, \label{Pbc1}\\
\Phi = 0 & \quad {\rm at} \quad r = 2, \label{Pbc2}
\end{eqnarray}
where ${\bm n}=(x,y)$ is the outward normal vector of the inner cylinder surface.

The no-slip conditions and Eq. (\ref{Pbc2}), which are Dirichlet boundary conditions,
are modeled by the classical VP method \cite{Angot1999}.
In contrast, 
Eq. (\ref{Pbc1}), which is an inhomogeneous Neumann boundary condition, 
can be modeled by the flux-based VP method we have proposed.
We here consider flow motions in the extended square domain $\Omega=\{(x,y)|-2.56 \leq x < 2.56, -2.56 \leq y < 2.56 \}$ 
under the periodic conditions at $x = \pm 2.56$ and $y= \pm 2.56$.
The resulting penalized equations read 
\begin{eqnarray}
\partial_t {\bm u} &=& -\nabla p - ({\bm u}  {\bm \cdot} \nabla) {\bm u}  + {\mathrm {Pr}} \nabla^2 {\bm u} + \left( 1-\chi \right){\mathrm {Ra Pr}} \Phi {\bm e}_y -\frac{\chi {\bm u}}{\eta_d}  , \label{eq5-2-1} \\
\partial_t \Phi &=& -(1-\chi)({\bm u} {\bm \cdot} \nabla) \Phi \nonumber \\
                &+& \nabla {\bm \cdot} \left[ \left\{(1-\chi_2)+\eta_n \chi_2 \right\} \nabla \Phi + \chi_2 {\bm \beta} \right] - \chi_2 \left(  \nabla {\bm \cdot}  {\bm \beta} \right) - \frac{\chi_1 \Phi}{\eta_d}, \label{eq5-2-3}
\end{eqnarray}
with Eq. (\ref{eq5-1-2}),
where $\eta_d \, (>0)$ and $\eta_n \, (>0)$ are the penalization parameters,
${\bm \beta}({\bm x})={\bm x}/r$ with ${\bm \beta}({\bm 0})={\bm 0}$,
and the mask function $\chi({\bm x})$ is defined by $\chi({\bm x}) = \chi_1({\bm x}) + \chi_2({\bm x})$ 
in which $\chi_1$ and $\chi_2$ are given by
\begin{equation}
\chi_1({\bm x}) = \left \{
\begin{array}{ll}
1   & {\rm if} \quad r > 2, \\ 
1/2 & {\rm if} \quad r = 2, \\
0   & {\rm otherwise},
\end{array}
\right.
\,\,\, {\mathrm {and}} \,\,\,
\chi_2({\bm x}) = \left \{
\begin{array}{ll}
1   & {\rm if} \quad r < 1, \\ 
1/2 & {\rm if} \quad r = 1, \\
0   & {\rm otherwise}.
\end{array}
\right. 
\label{eq5-4}
\end{equation}
We can select another expression for ${\bm \beta}$ as far as ${\bm n} {\bm \cdot} {\bm \beta} =1$ at $r=1$.
Note that $\nabla {\bm \cdot} (\chi_2 {\bm \beta}) - \chi_2 (\nabla {\bm \cdot} {\bm \beta}) =0 $ for $r <1$.
In Appendix \ref{appeC}, we discuss the VP representation of mixed Dirichlet--Neumann boundary conditions for a 1D Poisson equation.

We numerically solve Eqs. (\ref{eq5-2-1}) and (\ref{eq5-2-3}) with (\ref{eq5-1-2}), 
using the marker-and-cell method with a staggered grid at $256 \times 256$ grid points.
Second order central finite-differences and interpolation in Eqs. (\ref{2ndfd}) and (\ref{2ndint}) are employed for the spatial discretization, 
and the explicit Euler method is used as the time marching. 
The time increment $\Delta t$ is $10^{-6}$.
The computations were carried out at ${\mathrm{Pr}}=0.7$ for two ${\mathrm {Ra}}$ numbers, ${\mathrm {Ra}}=5700$ and ${\mathrm {Ra}}=5 \times10^4$.
The penalization parameters are taken as $\eta_d=\eta_n=5 \times 10^{-6}$.
The Poisson equation for the pressure $p$ is solved by the use of a successive over relaxation method. 
The convergence of $p$ is here achieved, 
if the relative residual error of $p$ measured by the $\ell^1$-norm becomes less than $10^{-6}$.
The numerical solutions are regarded as steady ones, 
if the relative time increments of $u$, $v$ and $\Phi$ normalized by $\Delta t$ respectively become less than $10^{-6}$,
where the increments are measured by the $\ell^1$-norm, and ${\bm u}=(u,v)$.

We also performed numerical simulations at ${\mathrm{Ra}} = 5 \times10^4$ and ${\mathrm{Pr}}=0.7$,
using different number of grid points, penalization parameters $\eta_n=\eta_d$, and time increment $\Delta t$.
The simulations confirmed that the numerical results for $128^2$ with $\eta_d=\eta_n=5 \times 10^{-6}$ and $\Delta t=10^{-6}$
excellently agree with those at $256^2$ using the same values of $\eta_d, \eta_n$, and $\Delta t$.
Such an excellent agreement of the numerical results is also observed, when doubling either $\Delta t$ or $\eta_d(=\eta_n)$. 
This convergence study justifies the choice of the parameters in the above computations.
Note that for stability reasons due to the explicit time discretization, 
we have chosen the values of $\Delta t$ so as to satisfy $\Delta t < \eta_d$, $\Delta t < \eta_n$ and $\Delta t < 1/{\mathrm{Ra}}$.

\subsection{Numerical results} \label{sec5-2}
Our numerical results are compared with results in Refs. \cite{Ren2013,Yoo2003}
in which also second order finite-differences are used.
Yoo \cite{Yoo2003} computed the flows in the polar coordinate system,
while Ren et al. \cite{Ren2013} used a Cartesian coordinate system combined with an immersed boundary method for Neumann boundary conditions which they developed.
Ren et al. \cite{Ren2013} reported that their results are in accordance with those of Ref. \cite{Yoo2003}.

\begin{figure}[tb]
\begin{center}
\includegraphics[width=7cm,keepaspectratio]{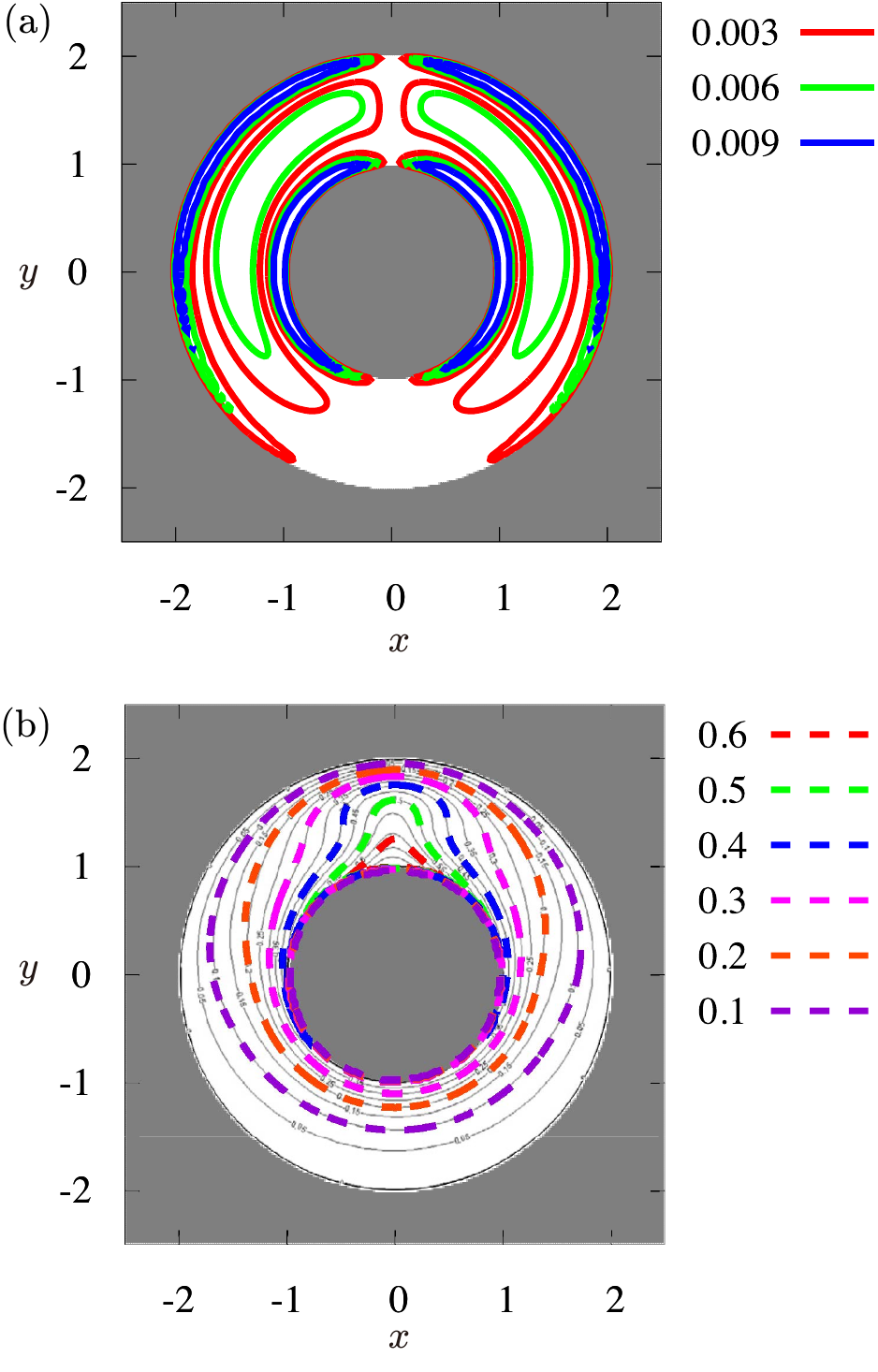}
\end{center}
\caption{(a) The streamlines and (b) the isotherms of $\Phi$ at ${\mathrm {Ra}}=5700$.
The isotherms denoted by the dashed lines are compared with those of Ref. \cite{Ren2013} that are denoted by the solid lines. 
The gray regions belong to the solid region $\Omega \backslash {\bar \Omega}_f$.} \label{fig5-1}
\end{figure}

Fig. \ref{fig5-1} shows the streamlines and the isotherms of $\Phi$ at ${\mathrm {R}}=5700$.
The streamlines and isotherms are symmetric with respect to $x=0$.
In Fig. \ref{fig5-1} (b), 
we observe that the isotherms at $0.1, \ 0.2, \ 0.3, \ 0.4, \ 0.5 $ and $0.6$ excellently agree with those presented in Fig. 7 (b) of Ref. \cite{Ren2013}.
We omitted the other isotherms of our results for brevity, since they also show excellent agreement with those of Ref. \cite{Ren2013}.

\begin{figure}[tb]
\begin{center}
\includegraphics[width=8cm,keepaspectratio]{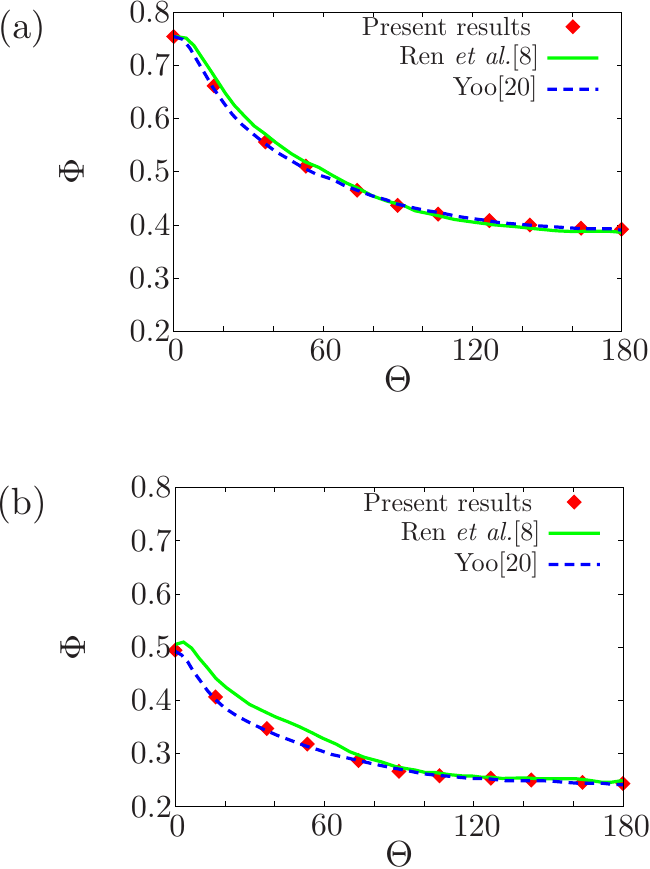}
\end{center}
\caption{Comparison of $\Phi$ on the inner cylinder surface at (a) ${\mathrm Ra}=5700$ and (b) ${\mathrm Ra}=5 \times 10^4$.}
\label{fig5-2}
\end{figure}

The spatial distributions of $\Phi$ on the inner cylinder at ${\rm Ra}=5700$ and $5 \times 10^{4}$ are shown in Fig. \ref{fig5-2},
 and are compared with the results in Refs. \cite{Ren2013,Yoo2003}. 
The angle $\Theta$ is measured counterclockwise from the top of the inner cylinder such that $\Theta=0^{\circ}$ at $(x,y)=(0,1)$ and $\Theta=180^{\circ}$ at $(x,y)=(0,-1)$. 
It can be seen that our results are in good agreement with those of Ref. \cite{Yoo2003} at both ${\rm Ra}$ numbers and the results of Ref. \cite{Ren2013} at ${\rm Ra}=5700$.
At ${\rm Ra}=5 \times 10^{4}$,
our result shows even better agreement with Ref. \cite{Yoo2003} than with Ref. \cite{Ren2013}.
The reciprocal of $\Phi$ at $r=1$ gives an average Nusselt number that characterizes the average heat transfer ratio through the inner cylinder surface \cite{Ren2013}.

{\it Note added in proof}: Guo et al. \cite{new} have developed an immersed boundary method for inhomogeneous Neumman boundary conditions in the frameworks of a finite volume method.
They used the same test case for validation, and showed that their numerical results likewise agree with those in Ref. \cite{Yoo2003}.


\section{Conclusions}\label{sec6}
We have developed a VP representation for inhomogeneous Neumann boundary conditions,
generalizing the technique proposed by Kadoch et al. \cite{KaKoAnSc} for homogeneous Neumann boundary conditions. 
It is based on the flux continuity between the fluid and solid domains.
The method was first applied to a 1D Poisson equation with inhomogeneous Neumann boundary conditions.
The analytical solution of the penalized Poisson equation converges to the exact solution of the non-penalized 1D Poisson equation for the penalization parameter $\eta \rightarrow 0$.
The numerical solutions of the penalized equation were obtained by the use of second order central finite-differences and interpolation.
In the case that the same values of the inhomogeneous Neumann boundary conditions are imposed at two interfaces between fluid and solid domains,
the errors between the exact solution of the non-penalized equation and the numerical solutions of the discretized equation exhibit second order convergence in terms of the grid width $h$ for sufficiently small $\eta$.
It was suggested that the $O(h^2)$ convergence is due to the cancellation of the discretized errors of the interface between the solid and fluid domains.
The eigenvalues of the discretized Laplace operator are identical for homogeneous and inhomogeneous boundary conditions and were studied in Ref. \cite{KoNgSc}.
We also showed that the VP representation needs a source term in the solid domain for the Poisson equation imposing two different values of inhomogeneous Neumann boundary conditions.
This term allows for the convergence of the penalized solution to the non-penalized one in the limit of $\eta \rightarrow 0$.
The numerical solutions of the penalized equation exhibit $O(h)$ convergence.
Our results suggest that the value of $\eta$ in the range of $10^{-5} \lesssim \eta \lesssim 10^{-3}$ is optimal for $h \approx 10^{-2}$.

The developed VP method was applied to the advection-diffusion equation coupled with the Navier--Stokes equations.
We simulated steady free convection of an incompressible fluid in a concentric annulus heated through the inner cylinder surface.
It was found that the VP method well preserves the temperature distributions on the inner cylinder  obtained in the literature \cite{Ren2013,Yoo2003}.

This work motivates us to develop a VP method further for Robin boundary conditions,
and apply the developed method to multiphysics problems in complex geometries, 
e.g., Stefan problems in crystal growth and to magnetohydrodynamic problems in magnetically confined fusion.

\section*{Acknowledgments}
The authors thank Dr. Dmitry Kolomenskiy for fruitful discussion about the convergence issue in the penalized Poisson equation with homogeneous Neumann boundary condition.
This work was partially supported by JSPS KAKENHI Grant Numbers (S)16H06339 and (C)17K05139.
K.S. acknowledges support by the French Federation for Magnetic Fusion Studies (FR-FCM) and of the Eurofusion consortium from the Euratom research and training programme 2014-2018 and 2019-2020 under grant agreement No 633053. The views and opinions expressed herein do not necessarily reflect those of the European Commission.

\clearpage
\appendix
\section{On different values of inhomogeneous Neumann conditions} \label{appeA}
To get deeper insight into the discretization error of the VP representation for inhomogeneous Neumann boundary conditions, 
we consider the Poisson equation (\ref{eq2-0}) with the source term $f(x)$ given by
\begin{eqnarray}
f(x) = (1-\epsilon) \cos{x} +\epsilon \sin{x},  \label{eqA-1}
\end{eqnarray}
in the fluid domain $\Omega_f =\{x \ | \ 0 \le x < \pi \}$ imposing the inhomogeneous Neumann boundary conditions
\begin{equation}
d_x w(0)=\alpha+\epsilon , \quad {\mbox{and}} \quad d_x w(\pi)=\alpha -\epsilon, \label{eqA-2}
\end{equation}
where $\epsilon$ and $\alpha$ are constant. 
Under the condition $\int_{\Omega_f} w (x) d x=0$, 
the exact solution is given by 
\begin{eqnarray}
w(x) = (1-\epsilon)\cos{x}+\epsilon \sin{x} + \alpha x - \frac{\pi \alpha}{2} - \frac{2 \epsilon}{\pi}. \label{eqA-3} 
\end{eqnarray}
Applying the VP representation to Eqs. (\ref{eqA-1}) and (\ref{eqA-2}),
we obtain
\begin{equation}
-d_x ( \theta d_x v + \chi \beta)= (1-\chi) f -\chi d_x \beta, \label{eqA-4}
\end{equation}
where $\beta(x)=\alpha + \epsilon \cos(x)$, $\Omega_f$ is extended to the domain $\Omega =\{x \ | \ 0 \le x <2 \pi \}$,
and then $2\pi$-periodic boundary conditions are imposed. 
The functions $\theta(x)$ and $\chi (x)$ are respectively given by Eqs. (\ref{eq2-5}) and (\ref{eq2-7}).
The analytical solution of Eq. (\ref{eqA-4}) is obtained in the same way as in Section \ref{sec2-1}, and is omitted here for brevity.

\begin{figure}[tb]
\begin{center}
\includegraphics[width=7cm,keepaspectratio]{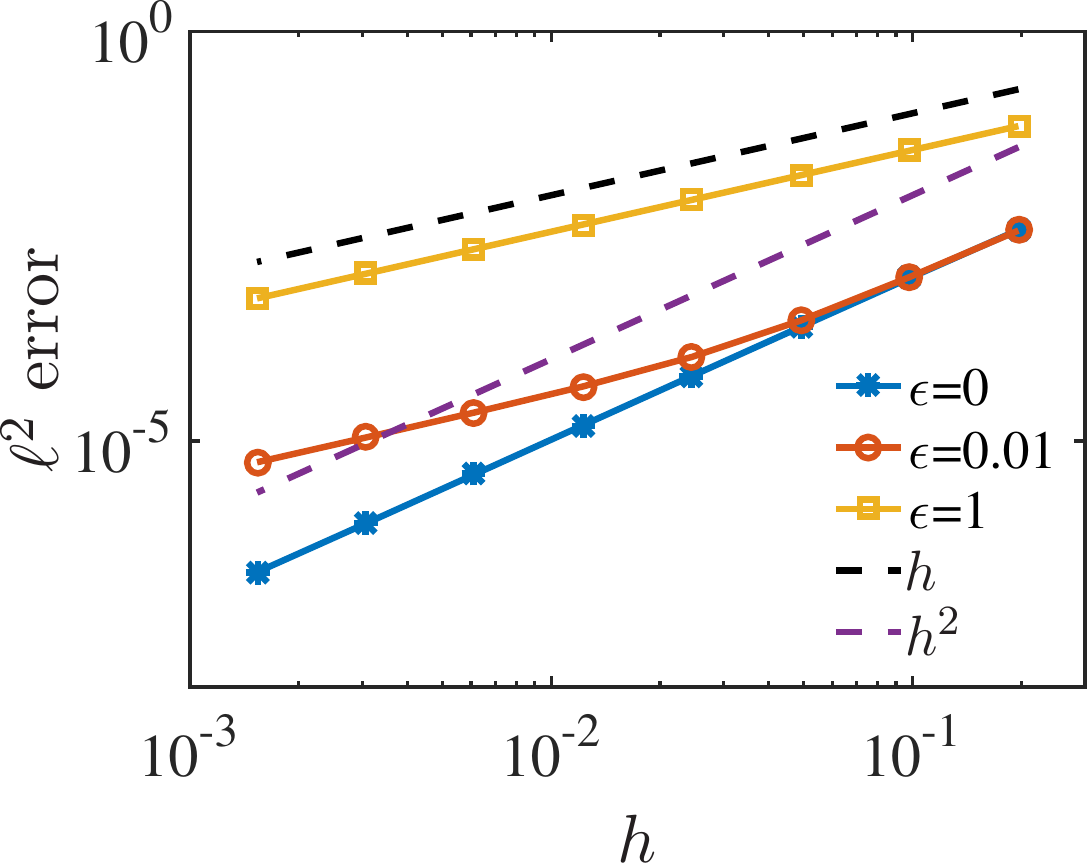}
\end{center}
\caption{The $h$-dependence of the $\ell^2$ error between the numerical solutions of Eq. (\ref{eqA-4}) and the exact solution (\ref{eqA-3}) for $\alpha=1$ and $\eta=10^{-8}$ using $\epsilon=0,\ 0.01$ and $1$.
The errors are computed in the fluid domain $\Omega_f$.}
\label{figA-2}
\end{figure}

Fig. \ref{figA-2} shows that the $\ell^2$ error between the numerical solutions $v$ of Eq. (\ref{eqA-4}) and the exact one (\ref{eqA-3}) for $\eta=10^{-8}$ using different $\epsilon$.
We see the $O(h^2)$ convergence of the error at $\epsilon=0$ which degrades to an $O(h)$ convergence for $\epsilon=1$.
It is observed at $\epsilon \,=0.01$ that this error decays approximately as $O(h^2)$ for large $h \, (\gtrsim 5 \times 10^{-2})$,
and they show the $O(h)$ decay for smaller $h \, (\lesssim 10^{-2})$.
The errors measured in the $\ell^1$- and $\ell^\infty$-norms show the similar behavior, 
and their figures are omitted for brevity.

\section{Penalized 1D thermal diffusion equation} \label{appB}
We here consider a time-dependent problem corresponding to a 1D heat equation in the domain $\Omega_f=\{x \ | \  -\pi < x < \pi \}$:
\begin{equation}
\partial_t \varphi (x,t) = \partial_x^2 \varphi(x,t), \label{eq4-1}
\end{equation}
with the initial condition given by 
\begin{eqnarray}
\varphi(x,0) = \sin{x}, \label{eq4-2-1} 
\end{eqnarray}
and the inhomogeneous time-dependent boundary conditions given by
\begin{eqnarray}
\partial_x \varphi(\pi ,t) = \partial_x \varphi(-\pi,t) = -e^{- t}, \label{eq4-2-2} 
\end{eqnarray}
where $\partial_t= \partial/\partial t$.
The exact solution then reads
\begin{eqnarray}
\varphi(x,t) = e^{- t} \sin x. \label{eq4-3}
\end{eqnarray}

The VP representation of Eqs. (\ref{eq4-1}) and (\ref{eq4-2-2}) results in
\begin{eqnarray}
\partial_t \varphi_p =\partial_x \left[ \left\{ (1-\chi) + \eta \chi \right\} \partial_x \varphi_p + \alpha \chi \right], \label{eq4-4}
\end{eqnarray}
where $\alpha = -e^{- t}$, and $\chi (x)$ is a mask function defined by 
\begin{eqnarray}
\chi(x)= \left\{
\begin{array}{ll}
0 & {\rm if} \quad -\pi < x < \pi ,  \\ 
1/2 & {\rm if} \quad x = \pm \pi , \\
1 & {\rm otherwise}  .
\end{array}
\right. 
\end{eqnarray}
Here, we employ a periodic domain $\Omega = \{x \ | \ -\pi-0.2 \leq x < \pi +0.2 \}$.
Based on Eqs. (\ref{eq4-2-2}) and (\ref{eq4-3}), 
we use the initial condition $\varphi_p(x,0) = ( 1-\chi ) \sin{x}$. 

\begin{figure}[!t]
\begin{center}
\includegraphics[width=6.5cm,keepaspectratio]{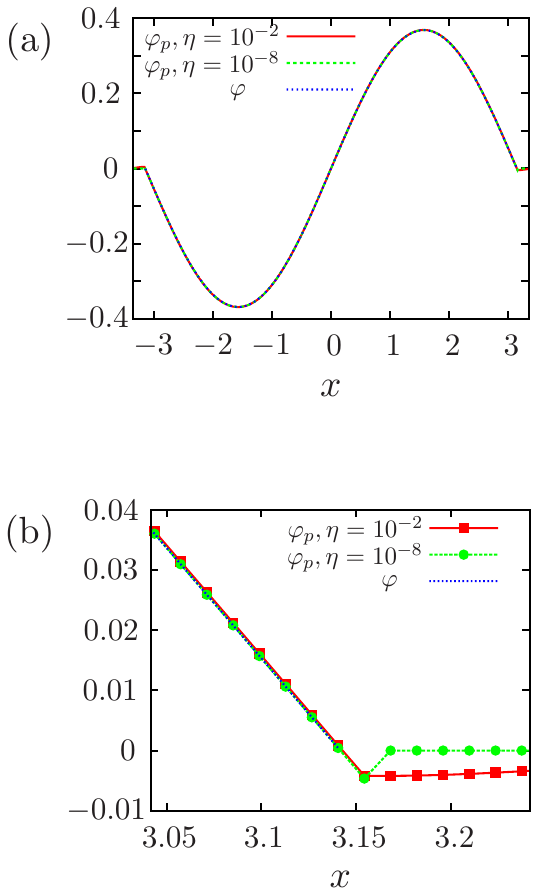}
\end{center}
\caption{(a) The numerical solution $\varphi_p(x,t)$ of Eq. (\ref{eq4-4}) for $\eta=10^{-2}$ and $10^{-8}$ and the exact solution $\varphi(x,t)$, Eq. (\ref{eq4-3}), at $t=1$.
(b) A zoom of Fig. \ref{fig4-1} (a) around $x=\pi$.}
\label{fig4-1}
\end{figure}

\begin{figure}[!t]
\begin{center}
\includegraphics[width=7cm,keepaspectratio]{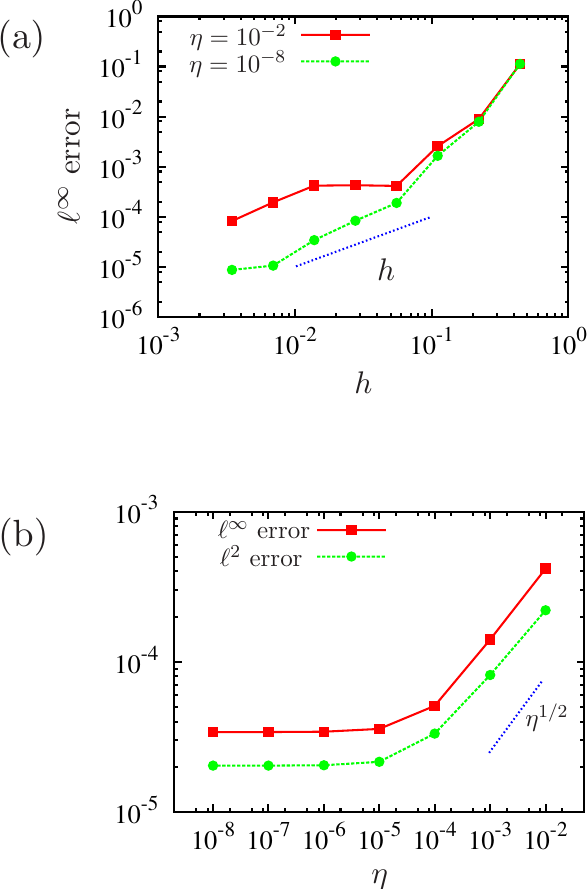}
\end{center}
\caption{(a) The $h$-dependence of the $\ell^{\infty}$ errors between the numerical solutions $\varphi_p$ and the exact solution $\varphi$, Eq. (\ref{eq4-3}), at $t=1$ for $\eta=10^{-2}$ and $10^{-8}$.
(b) The $\eta$-dependence of the $\ell^{\infty}$ and $\ell^{2}$ errors between the solutions $\varphi_p$ and $\varphi$ at $t=1$.
The errors are computed in the fluid domain $\Omega_f$.}
\label{fig4-2}
\end{figure}

Equation (\ref{eq4-4}) was numerically solved, using second order central finite-differences and interpolation with $N=512$ grid points.
The Crank-Nicolson method is employed for time integration with time increment $\Delta t = 10^{-5}$. 
The penalization parameter $\eta$ is set to be either $10^{-2}$ or $10^{-8}$. 
We obtain the numerical solution $\varphi_p$, 
which is compared with the exact solution $\varphi(x,t)$ of Eq. (\ref{eq4-3}).
Fig. \ref{fig4-1} shows $\varphi_p$ and $\varphi$ at $t=1$. 
In $\Omega_f$, $\varphi_p$ excellently agrees with $\varphi$.
To quantify the precision, 
we present the $\ell^{\infty}$ and $\ell^{2}$ errors between $\varphi_p$ and $\varphi$ at $t=1$ in Fig. \ref{fig4-2}.
In Fig. \ref{fig4-2} (a), it can be seen that the errors decay, as the grid width $h$ decreases.
Here, in this section, $h=(2\pi+0.4)/N$.
For smaller $h$, we observe a $O(h)$ decay, since no grid point is on the boundary of $\Omega_f$.
Fig. \ref{fig4-2} (b) depicts the $\eta$-dependence of the errors. 
The errors decay approximately as $O(\eta^{1/2})$ for larger $\eta \, (\gtrsim 10^{-4})$ and converge for smaller $\eta \, (\lesssim 10^{-5})$.
The decay of $O(\eta^{1/2})$ has been proved in theorem 2.1 of Ref. \cite{KaKoAnSc}.


\section{Mixed Dirichlet--Neumann boundary conditions for the 1D Poisson equation} \label{appeC}
We consider here the 1D Poisson equation (\ref{eq2-0}) with the source term given by
\begin{eqnarray}
f(x) = \cos{x},  \label{eqC-1}
\end{eqnarray}
in the domain $\Omega_f=\{x \ | \ 0 < x < \pi\}$ imposing mixed Dirichlet and inhomogeneous Neumann boundary conditions, i.e., 
\begin{equation}
w(0)=0 , \quad {\mbox{and}} \quad  d_x w(\pi)=\alpha, \label{eqC-2}
\end{equation}
with $\alpha$ being constant.
The exact solution of Eq. (\ref{eqC-1}) reads
\begin{eqnarray}
w(x) = \cos{x} + \alpha x -1. \label{eqC-3}
\end{eqnarray}

Applying the VP representation to Eqs. (\ref{eqC-1}) and (\ref{eqC-2}),
we obtain
\begin{equation}
-d_x ( \theta d_x v + \chi_n \beta)= (1-\chi) f -\frac{\chi_d v}{\eta_d} -\chi_n d_x \beta ,
\label{eqC-6}
\end{equation}
in which $\beta (x) = \alpha \left( 1+ \sin x \right)$, 
$\theta = \left( 1- \chi_n \right) + \eta_n \chi_n$,
and $\chi(x)=\chi_d(x) + \chi_n(x)$ where
\begin{eqnarray}
& & \chi_d(x) = \left\{
\begin{array}{ll}
0   & {\rm if} \quad 0 \le x < 3\pi/2, \\
1/2 & {\rm if} \quad x=0, 3\pi/2, \\
1   & {\rm otherwise} ,  \\
\end{array}
\right.
\end{eqnarray}
and
\begin{eqnarray}
& & \chi_n(x) = \left\{
\begin{array}{ll}
0   & {\rm if} \quad 0\le x < \pi, \,\,  3\pi/2 <x <2\pi, \\
1/2 & {\rm if} \quad x=\pi, 3\pi/2, \\
1   & {\rm otherwise} .  \\
\end{array}
\right. 
\label{eqC-7}
\end{eqnarray}
The domain $\Omega_f$ is extended to $\Omega=\{x \ | \ 0 \le x <2\pi \}$,
and $2\pi$-periodic boundary conditions are imposed.

The analytical solution of Eq. (\ref{eqC-6}) is uniquely obtained not only by using 
$v(\pi^{-})=v(\pi^{+})$, $v(3\pi/2^{-})=v(3\pi/2^{+})$, $v(0^{+})=v(2\pi^{-})$,
and the continuity of the flux, $d_x v = \eta d_x v + \beta$, at $x= \pi$ and $3\pi/2$,
but also by imposing $d_x v (0^+)=d_x v (2\pi^{-})$.
For brevity we omit the form of the analytical solution, because it is a rather complicated expression. 
The analytical solution converges to the exact solution (\ref{eqC-3}) for $\eta_d \rightarrow 0$ and $\eta_n \rightarrow 0$.

\begin{figure}[!t]
\begin{center}
\includegraphics[width=7cm,keepaspectratio]{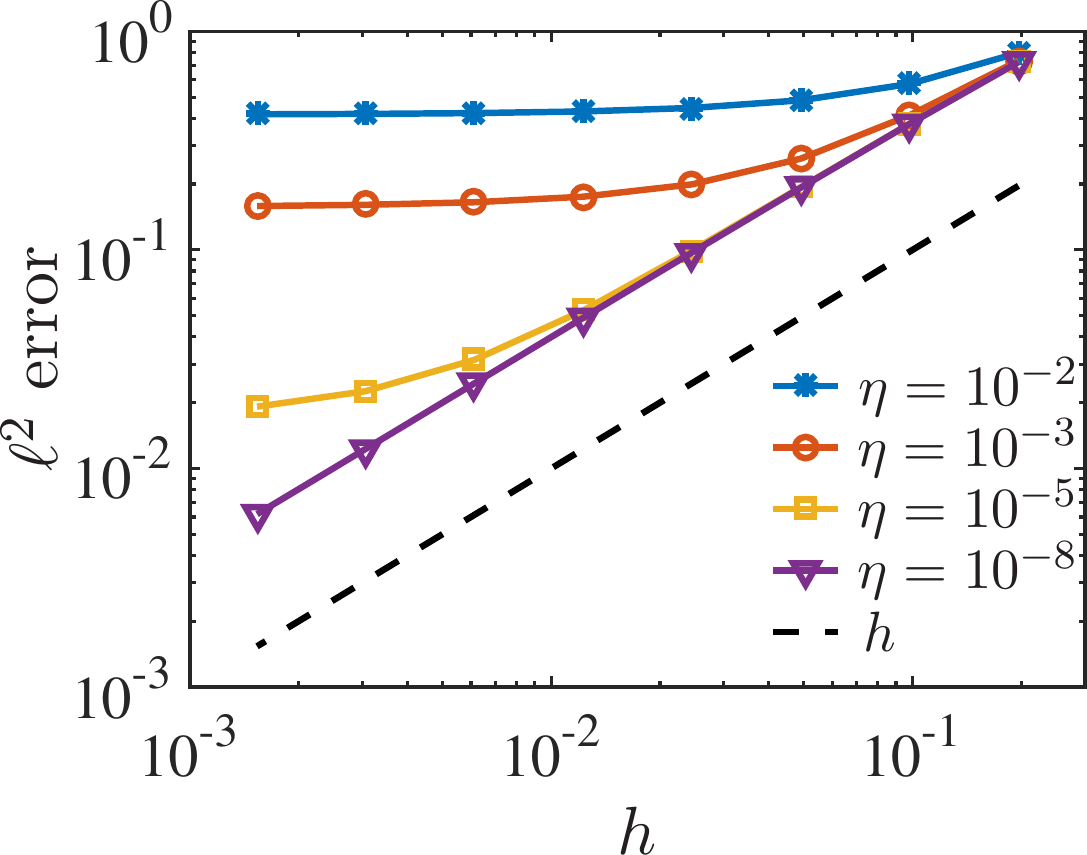}
\end{center}
\caption{The $h$-dependence of the $\ell^2$ errors between the numerical solutions of Eq. (\ref{eqC-6}) and the exact solution of Eq. (\ref{eqC-1}) for $\alpha=1$ using $\eta=10^{-2},\ 10^{-3},\ 10^{-5}$, and $10^{-8}$. The errors are computed in the fluid domain.} \label{figC-2}
\end{figure}

Fig. \ref{figC-2} shows that the errors between the numerical solutions of Eq. (\ref{eqC-6}) and the exact solution $w$, Eq. (\ref{eqC-3}), for $\alpha=1$ using different $\eta_d =\eta_n= \eta$.
For $\eta =10^{-3}$, the errors decay approximately as $O(h)$ with decreasing $h$, and then saturate due to the penalization error.
The level of the saturation becomes smaller, as $\eta$ decreases.
It can be seen that for sufficiently small $\eta$ the errors decay like $O(h)$.

\end{document}